\documentclass[12pt,a4paper,twoside,final,notitlepage, leqno]{article}
\usepackage[english]{babel}
\usepackage[T1]{fontenc}  
\usepackage{epsfig, graphicx, amssymb}
\usepackage{amsmath,amsthm,epsfig,amsfonts}  
\usepackage{float}
\usepackage{color}
\def\br {\break}

\def \smb {{\scriptstyle \bullet }}
\newcommand{\monitem}{ \smallskip \noindent $\bullet$ \quad  } 
\newcommand{\moneq}{\vspace*{-7pt} \begin{equation} \displaystyle } 
\newcommand{\moneqstar}{\vspace*{-6pt} \begin{equation*} \displaystyle } 
\newcommand{\monendstar}{\vspace*{-6pt} \end{equation*}   }
\newcommand{\monend}{\vspace*{-7pt} \end{equation}   }
\newcommand{\moneqarraystar}{ \begin{eqnarray*} \displaystyle } 
\newcommand{\monendarraystar}{ \end{eqnarray*}   }

\newcommand{\RR}[0]{\mathbb{R}}

\def\v{\vert}
\def\eq{\mathop{\rm {eq}}\nolimits}

\def\section*#1{}

\usepackage{fancyhdr}
\renewcommand{\headrulewidth}{0pt}

\begin{document} 

\fancypagestyle{plain}{ \fancyfoot{} \renewcommand{\footrulewidth}{0pt}}
\fancypagestyle{plain}{ \fancyhead{} \renewcommand{\headrulewidth}{0pt}}

~

  \vskip 2.1 cm

\centerline {\bf \LARGE  On anti bounce back boundary condition  }

\bigskip

\centerline {\bf \LARGE   for lattice Boltzmann schemes  }

 \bigskip  \bigskip \bigskip

\centerline { \large    Fran\c{c}ois Dubois$^{ab}$, Pierre Lallemand$^{c}$,  
 Mohamed-Mahdi Tekitek$^{d}$}

\smallskip  \bigskip 

\centerline { \it  \small   
$^a$   Dpt. of Mathematics, University Paris-Sud,  B\^at. 425, F-91405  Orsay, France.} 

\centerline { \it  \small   
$^b$    Conservatoire National des Arts et M\'etiers, LMSSC laboratory,  F-75003 Paris, France.} 

\centerline { \it  \small  $^c$   Beijing Computational Science Research Center, 
Haidian District, Beijing 100094,  China.}

\centerline { \it  \small  $^d$   Dpt. Mathematics, Faculty of Sciences of Tunis, 
University Tunis El Manar, Tunis, Tunisia. } 


\bigskip  \bigskip  

\centerline { 05 April 2019  
{\footnote {\rm  \small $\,$   Contribution  presented to the 
14th ICMMES  Conference, Nantes (France), 18 - 21 July 2017. 
{\it    Computers And Mathematics With Applications},
volume 79, pages 555-575, February 2020.}}}



\bigskip  \bigskip   
\noindent {\bf \large Abstract} 

\noindent 
In this contribution, we recall the derivation of the  anti bounce back boundary condition
for the D2Q9 lattice Boltzmann scheme. 
We recall various elements of the state of the art for anti bounce back applied to linear heat  and acoustics equations
and in particular the possibility to take into account curved boundaries.
We present an asymptotic analysis that allows an expansion of all the fields in the boundary cells.
This analysis based on the Taylor expansion method  confirms the well known behaviour  of  anti bounce back boundary for the heat equation.
The analysis puts also in evidence a hidden differential boundary condition in the case of linear acoustics.
Indeed, we observe discrepancies 
in the first layers near the boundary.
To reduce these discrepancies, we propose a new boundary condition mixing bounce back for the oblique links and
anti bounce back for the normal link. This  boundary condition is able to  enforce both pressure
and tangential velocity on the boundary. Numerical tests for the Poiseuille flow illustrate our theoretical analysis
and show improvements in the quality of the flow. 

\bigskip \bigskip   \newpage \noindent {\bf \large    1) \quad  Introduction  }    

\fancyhead[EC]{\sc{ Fran\c{c}ois Dubois, Pierre Lallemand, Mohamed-Mahdi Tekitek }} 
\fancyhead[OC]{\sc{On anti bounce back boundary condition for lattice Boltzmann schemes}} 
\fancyfoot[C]{\oldstylenums{\thepage}}

\noindent 
From the early days of Lattice Boltzmann 
scheme studies, boundary conditions have been the object of various proposals. 
In fact, two popular methods of boundary conditions are used to impose  given macroscopic conditions 
(velocity, pressure, thermal fields, \dots). The first one, called ``half-way'' 
was proposed in the context of cellular automata \cite {dHPL85}.   
The  ``half-way'' approach \cite{dHG03}, consists in using 
the distributions that leave the domain to define the unknown distributions by a simple reflection 
(called bounce back) 
or antireflection (called anti bounce back). 
In this method the physical position of the wall is located between the last internal domain
node and first node beyond boundary. In \cite{GA94,DLT10} the exact position 
of the physical wall is investigated for some simple problems.    
The second method proposed first by Zou and He \cite{ZH97}, called ``boundary on the node'', 
uses the projection of the given macroscopic conditions in the distribution 
and assumes the bounce back rule for the non-equilibrium part of the particle distribution.
 In this method the physical wall is located 
on the last node of the domain. 

\noindent 
Here, we only focus on ``half-way'' boundary conditions method. In a previous work~\cite{DLT15} a novel method
of analysis, based on Taylor developments, is used for  the bounce back scheme.
This linear analysis gives an expansion of the macroscopic quantity, on the physical wall, as powers of the mesh 
size. Later in~\cite{DLT17}  a generalized bounce back boundary condition for 
the nine velocities two-dimensional (D2Q9) lattice Boltzmann scheme is proposed.
This scheme is exact up to second order by elimination of spurious density terms (first order terms).

\noindent 
In this contribution, the anti bounce back boundary lattice Boltzmann scheme is investigated. 
This scheme is used to impose  Dirichlet boundary conditions for
the ``thermal problem'', {\it id est} the heat equation  where one scalar moment is conserved,  
or to impose a given pressure for a linear fluid problem. Many works proposed new  
boundary conditions~\cite{AK02, CL07,  FRJ02, LCM08,  WP02} intended to yield improved 
accuracy compared to the anti bounce back boundary scheme. 

\noindent   
In this paper,  the D2Q9 lattice Boltzmann scheme is first briefly introduced for heat equation and acoustic system  (Section 2).  
Then in Section 3, the anti bounce back boundary condition 
is presented to impose a given thermal field for thermal problem or  
to impose a given density/pressure for the linearized fluid. 
In Section 4, the scheme is analyzed using Taylor expansion for the heat equation. 
This asymptotic analysis gives an expansion of the conserved 
scalar moment which is exact up to order one.
In Section~5, we present the extended  anti bounce back that allows to handle curved boundaries. 
In Section 6, anti bounce back is analyzed for the linear fluid problem.
A hidden differential boundary condition is put in evidence. 
In Section 7, a mixing of bounce back and anti bounce back boundary scheme \cite{DL15}, is used and analyzed
to impose a given density and a given velocity on the physical wall.
Finally in Section 8 a novel boundary scheme is introduced and analyzed 
to impose a given pressure and tangential velocity on the boundary.    
A Poiseuille test case gives convergent results.

\newpage 
 \bigskip \bigskip   \noindent {\bf \large    2) \quad  
  D2Q9 lattice Boltzmann scheme for heat  and acoustic problems }   

\begin{figure}    [H]  \centering 
  \centerline { \includegraphics[width=.35 \textwidth, angle=0] {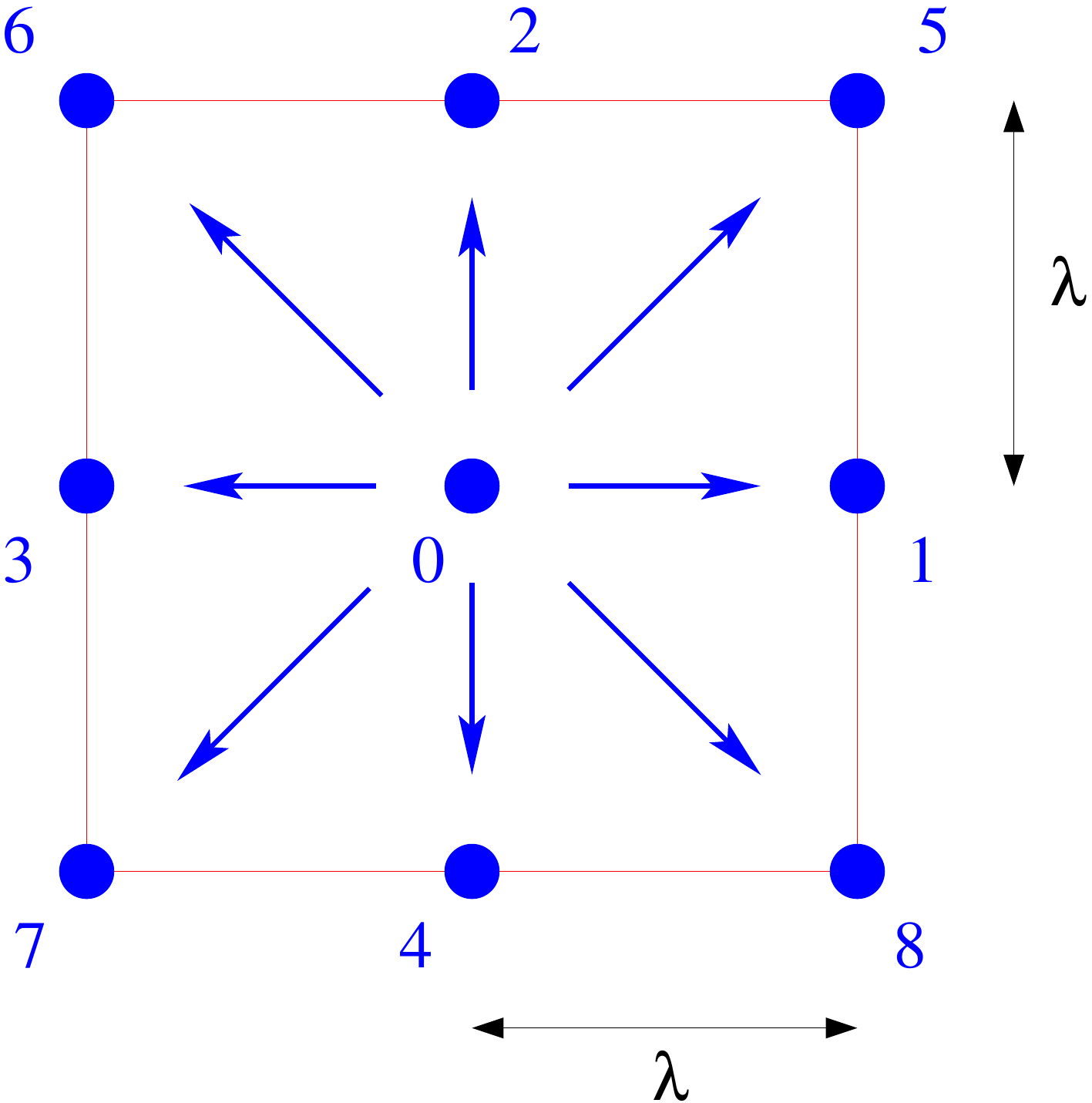}}
\caption{D2Q9 lattice Boltzmann scheme; discrete velocities $ \, v_j \, $  for $ \, 0 \leq j \leq 8 $ } 
\label{d2q9stencil} \end{figure}

\noindent  
The D2Q9 lattice Boltzmann scheme  uses a set of discrete velocities described 
in Figure~\ref{d2q9stencil}. 
A density distribution $ \, f_j \,$ is associated to each basic velocity~$ \, v_j  $.  
The two components and the modulus of this vector are denoted by $ \, v_{j x} $,  $ \, v_{j y} \, $ and $ \,  \v v_j \v \, $ respectively.   
The first three moments for the 
density and momentum are defined according to 
\moneq \label{moments-conserves} 
\rho \, = \,  \sum^8_{j=0} f_{j} \, = \, m_0 \,, \,\,\,  
J_x  \, = \,   \sum^8_{j=0}   v_{j x}  \, f_{j} \, = \, m_1 \,, \,\,\,    
J_y  \, = \,   \sum^8_{j=0} v_{j y}  \, f_{j} \, = \, m_2 \, .  
\monend    
We complete this set of moments and 
construct a vector  $ \, m \, $  of moments  $ \, m_k \,$ 
using Gram-Schmidt orthogonalization 
and an appropriate standardization to have simple expressions of the moment matrix, 
as proposed in \cite{LL00}:    
\moneq \label{momdents-ll00} 
  \left \{ \begin {array}{l} 
\displaystyle \varepsilon  =  3 \, \sum_{j=0}^8 \v v_j \v^2 \, f_j \,  
- 4 \,\, \lambda^2  \sum_{j=0}^8  \,  f_j \,, \\    \displaystyle 
\varphi_x  = \sum_{j=0}^8  \big[ (v_{j x} )^2  - (v_{j y}  )^2 \big]  \, f_j  \,, \,\, 
 \varphi_y  =   \sum_{j=0}^8   v_{j x}  \, v_{j y}    \,\,  f_j \,, \\ \displaystyle 
 q_x  =   3 \, \sum_{j=0}^8  \v v_j \v^2 \,  v_{j x} \, f_j 
 - 5 \, \lambda^2 \, \sum_{j=0}^8  \,   v_{j x} \, f_j \,,\,\,  
q_y   =   3 \, \sum_{j=0}^8  \v v_j \v^2 \,  v_{j y} \, f_j
-  5 \, \lambda^2 \, \sum_{j=0}^8  \, v_{j y} \,  f_j \,, \\\displaystyle 
  D = {9\over2} \, \sum_{j=0}^8  \v v_j \v^4 \, f_j
   \,\,  - \, {21\over2} \,  \lambda^2 \,  \sum_{j=0}^8  \v v_j \v^2   \,  f_j \,
  \,\, + \, 4 \,  \lambda^4 \, \sum_{j=0}^8  \,  f_j  \, . 
\end {array}  \right. \monend 
  %
%
The entire vector $ \, m \in \RR^9 \, $ of moments is defined by 
\moneq     \label{definition-moments-d2q9}
m  = \big( \rho \,,\, J_x \,,\, J_y ,\, \varepsilon \,,\, 
\varphi_x  \,,\, \varphi_y \,,\, q_x \,,\, q_y \,, \, D \big)^{\rm \displaystyle t} 
\monend 
and the previous relations can be written in a synthetic way~:   
\moneq \label{moments} 
m = M \, f \,. 
\monend
The invertible fixed  matrix $ \, M \,$ is usually \cite{LL00}  given by 
\moneq \label{matrice-M} 
 M  \, = \, \left (\begin {array}{ccccccccc}
\displaystyle   1 &   1 & 1 & 1 & 1 & 1 & 1 & 1 & 1\cr
\displaystyle   0 &   \lambda & 0 & -\lambda & 0 & \lambda & -\lambda & -\lambda  & \lambda\cr
\displaystyle   0 &   0 & \lambda & 0 & -\lambda & \lambda & \lambda & -\lambda & -\lambda\cr
\displaystyle -  4 \, \lambda^2 &   -\lambda^2  &  -\lambda^2 &  -\lambda^2 & -\lambda^2 &
 2 \,\lambda^2  &  2 \,\lambda^2 &  2 \,\lambda^2 &  2 \,\lambda^2 \cr  
\displaystyle   0 &    \lambda^2 & -  \lambda^2 &  \lambda^2 & -  \lambda^2 & 0 & 0 & 0 & 0\cr
\displaystyle   0 &   0 & 0 & 0 & 0 & 
\lambda^2 & -\lambda^2 & \lambda^2 & -\lambda^2    \cr  
\displaystyle   0 &   -2 \, \lambda^3 & 0 & 2  \, \lambda^3 & 0 & \lambda^3 &
 - \lambda^3 & - \lambda^3 &  \lambda^3 \cr
\displaystyle   0 &   0 & -2 \, \lambda^3 & 0 & 2 \, \lambda^3 &
  \lambda^3 &  \lambda^3 & - \lambda^3 & -  \lambda^3 \cr
\displaystyle   4 \, \lambda^4 &  -2 \, \lambda^4 &  -2  \, \lambda^4& -2  \, \lambda^4 &
 -2  \, \lambda^4 &  \lambda^4 & \lambda^4  &  \lambda^4 &   \lambda^4   \end{array}  \right) \, , 
\monend
where 
 $ \, \lambda \equiv {{\Delta x}\over{\Delta t}} \, $ is the  fixed numerical lattice velocity.
For scalar lattice Boltzmann  applications to thermal problems, the density $ \, \rho \, $   
is the   ``conserved variable''.

\monitem 
We suppose the following 
linear equilibrium of nonconserved moments 
 and associated relaxation coefficients. 
Due to (\ref{matrice-M}) and Table \ref{table-1}, we can explicit the vector $ \, f^{\eq} \,$ 
of  equilibrium values 
of the particle distribution~: 
\moneq \label{feq-d2q9-thermique} 
  \left \{ \begin {array}{l}  \displaystyle 
f_0^{\rm eq} =   {{\rho} \over{9}}    \, [ 1-\alpha +\beta ] \,,\,\,
f_1^{\rm eq} =   {{\rho} \over{36}}  \, [ 4 -\alpha-2\, \beta  ] \,,\,\,
f_2^{\rm eq} =   {{\rho} \over{36}}  \, [  4 -\alpha-2 \,\beta  ] \,, \\ \displaystyle  \vspace{-4 mm} \\ \displaystyle 
f_3^{\rm eq} =   {{\rho} \over{36}}  \, [ 4 -\alpha-2 \, \beta  ] \,,\,\,
f_4^{\rm eq} =   {{\rho} \over{36}} \,  [ 4 -\alpha-2 \, \beta ] \,,\,\, 
f_5^{\rm eq} =   {{\rho} \over{36}}  \, [ 4 +2 \, \alpha+\beta ] \,, \\ \displaystyle \vspace {-4mm}   \\ \displaystyle
f_6^{\rm eq} =   {{\rho} \over{36}}  \, [ 4 +2 \, \alpha+\beta   ] \,,\,\, 
f_7^{\rm eq} =   {{\rho} \over{36}}  \, [ 4 +2 \, \alpha+\beta    ] \,,\,\, 
f_8^{\rm eq} =   {{\rho} \over{36}} \, [ 4 +2 \, \alpha+\beta   ] \, . 
 \end{array} \right.  
\monend 
%
%
\begin{table}  [H]     \centering
 \centerline { \begin{tabular}{|c|c|c|c|c|c|c|c|c|}    \hline 
 moment &  $ J_x $ & $ J_y $ & $  \varepsilon $ &  $\varphi_x $ & $ \varphi_y $ & $ q_x $ & $  q_y $ & $ D $ \\   \hline 
 equilibrium  & $ 0 $ & $ 0 $ & $  \alpha  \, \lambda^2  \, \rho  $ & $0$ & $0$  & $0$ & $0$ & 
$  \beta \, \lambda^4 \, \rho $    \\   \hline 
 relaxation coefficient   & $ s_j $  & $ s_j $ &  $ s_e $ &  $ s_x $  &  $ s_x $ & $ s_q $ &   $ s_q $ 
&   $ s_d $ \\   \hline  
\end{tabular} }  
\caption{D2Q9 linear equilibrium of nonconserved moments 
 and associated relaxation coefficients for a thermal type  problem. } 
\label{table-1} \end{table}
%

\monitem 
During the relaxation step  of the heat diffusion problem, the conserved variable $ \, \rho \, $
is  not   modified.
In the framework of multiple relaxation times \cite{DDH92}, 
the nonconserved moments $\, m_1 \, $ to  $\, m_8 \, $ relax towards 
an equilibrium value  $ \, m_k^{\eq} \, $ displayed in Table \ref{table-1}. 
The relaxation step $ \,  m \longrightarrow m^*   \,$ needs also parameters
$ \, s_k \, $ presented also in the Table   \ref{table-1}: 
\moneq    \label{relaxation-moments}
  m_k^* \, = \,  m_k \,+\, s_k \, \big(  m_k^{\eq} \,-\,  m_k \big)  \,, \qquad k \geq 1 \, . 
\monend
Companion parameters   $ \, \sigma_k \,$  have been introduced by H\'enon  \cite {He87}   
in the context of cellular automata:  
\moneqstar 
 \sigma_k \, = \, {{1}\over{s_k}} \,-\, {1\over2} \, . 
\monendstar 
In terms of moments, the linear collision $ \, m \longrightarrow m^* \, $ 
can be written as  
\moneq \label{matrice-C}
m^* \,=\, C \,\, m \,,
\monend 
with a collision matrix $ \, C \, $ given by 
\moneq  \label{matrice-collision-thermic}
\!\!\! C \!=\! \left (
\begin{array}{ccccccccc}
\!\!\!\! 1 & \!\!\!\!\! 0 & \!\!\! 0 & \!\!\! 0 &  \!\!\! 0 &   \!\!\!\!\! 0 &  \!\!\!\!\! 0 &
\!\!\!\! 0\!\!\!\! &\!\!\!\! 0 \!\!\!\! \\ 
\!\!\!\!\!  0 & \!\!\!\!\! 1 \!\! - \!\! s_j  
&  \!\!\! 0 &  \!\!\! 0 &  \!\!\! 0 &   \!\!\!\!\! 0 &  \!\!\!\!\! 0 &
\!\!\!\! 0\!\!\!\! &\!\!\!\! 0 \!\!\!\! \\ 
\!\!\!\!\!  0 & \!\!\!\!\! 0 &  \!\!\!  1 \!\! - \!\! s_j  
 &  \!\!\! 0 &  \!\!\! 0 &   \!\!\!\!\! 0 &  \!\!\!\!\! 0 &
\!\!\!\! 0\!\!\!\! &\!\!\!\! 0 \!\!\!\! \\ 
\!\!\!\!\!  \alpha \, s_e \, \lambda^2 & \!\!\!\!\! 0 &  \!\!\! 0 &  \!\!\!\!\! 1-s_e &  \!\!\!\!\! 0 &  
 \!\!\! 0 &  \!\!\! 0 & \!\!\!\! 0 \!\!\!\!&\!\!\!\! 0 \!\!\!\! \\ 
\!\!\!\!\!  0 & \!\!\! 0 &  \!\!\!\! 0 &  \!\!\! 0 &  \!\!\! 1-s_x  &   \!\!\!\!\! 0 &  \!\!\!\!\! 0 
&\!\!\!\! 0 \!\!\!\!&\!\!\!\! 0 \!\!\!\! \\ 
\!\!\!\!\!  0 & \!\!\! 0 &  \!\!\!\! 0 &  \!\!\! 0 &  \!\!\! 0  &    \!\!\!\!\! 1-s_x &  \!\!\!\!\! 0 
&\!\!\!\! 0 \!\!\!\!&\!\!\!\! 0  \!\!\!\! \\ 
\!\!\!\!\!  0 & \!\!\!\!\! 0 
&  \!\!\! 0 &  \!\!\! 0 &  \!\!\!\!\! 0 &   \!\!\!\!\! 0 & 
\!\!\!\!1-s_q &   \!\!\! 0 \!\!\!\!&\!\!\!\! 0 \!\!\!\! \\ 
\!\!\!\!\!  0 & \!\!\!\!\! 0 &  \!\!\! 0 
&   \!\!\! 0 &  \!\!\!\!\! 0 &   \!\!\!\!\! 0 &  
\!\!\! 0 & \!\!\!\!1-s_q \!\!\!\!&\!\!\!\! 0 \!\!\!\! \\ 
 \beta \, s_d \, \lambda^4 & \!\!\!\!\! 0 &  \!\!\! 0 &  \!\!\! 0 &  \!\!\! 0 &   
\!\!\!\!\! 0 &  \!\!\!\!\! 0 &\!\!\!\! 0\!\!\!\! &  \!\!\!\! 1-s_d 
\end{array} \right) \, . \monend 
The iteration of the D$2$Q$9$ scheme  after the relaxation step follows the relation    
\moneqstar 
f^*_i(x, \, t) \, = \, \sum_\ell \, M^{-1} _{i \, \ell} \,  m_\ell^{*} \,, \,\, 0 \, \leq \, i \, \leq \, 8 \,  
\monendstar 
and  the  population  $ \, f_j (x,\, t+ \Delta t) \, $ ($ 0\leq j \leq 8$) 
at the new time step is  evaluated according to    
\moneq     \label{lb-iteration}
   f_j(x,t+\Delta t) \, = \, f^*_j(x - v_j\Delta t,t) \, , \,\,     0\leq j \leq 8 \, . 
\monend 

\monitem 
In \cite {Du08} and \cite{DL09}, we have introduced 
classical Taylor expansions in order to recover equivalent partial differential 
equations of the previous internal scheme. 
For the   thermal model, when $\, \Delta x \, $ and $\, \Delta t \, $ 
tend to zero, with a fixed ratio $ \, \lambda \equiv {{\Delta x}\over{\Delta t}} \, $
and fixed relaxation coefficients $ \, s_k $, 
 the conserved variable $\, \rho \, $ is solution of an asymptotic heat equation
\moneq    \label{chaleur}
{{\partial \rho}\over{\partial t }} - \, \mu \, \Delta \rho 
= {\rm O}(\Delta x^2)  \, . 
\monend 
The infinitesimal diffusivity  $ \,  \mu \, $ is evaluated according to 
\moneqstar 
\mu = {1\over6} \, \sigma_j \, (\alpha + 4 ) \,  \lambda \, \Delta x \, . 
\monendstar 

\monitem 
If we consider a linear fluid problem, the three first moments $ \, \rho$, $ \, J_x \equiv \rho \, u_x\,$
and $ \, J_y \equiv \rho \, u_y \, $ are conserved during the collision step. The table of equilibria and relaxation
 coefficients is modified in the  way specified in Table~\ref{table-2}:    
%
\begin{table}  [H]     \centering
 \centerline { \begin{tabular}{|c|c|c|c|c|c|c|c|c|}    \hline 
 moment & $  \varepsilon $ &  $\varphi_x $ & $ \varphi_y $ & $ q_x $ & $  q_y $ & $ D $ \\   \hline  
 equilibrium  & $  \alpha  \, \lambda^2  \, \rho  $ & $0$ & $0$ & 
$ - \lambda^2 \, J_x  $ & $ - \lambda^2 \, J_y $ &   $  \beta \, \lambda^4 \, \rho $    \\   \hline 
 relaxation coefficient    &  $ s_e $ &  $ s_x $  &  $ s_x $ & $ s_q $ &   $ s_q $ &   $ s_d $ \\   \hline  
\end{tabular} }  
\caption{
D2Q9 linear equilibrium of nonconserved moments 
 and associated relaxation coefficients for a acoustic linear fluid  problem. }   
\label{table-2} \end{table}
%
\vspace{-4 mm}
\noindent We can also display  for the linear fluid  case 
all components of  the linear vector $ \, f^{\eq} $:
\moneq \label{feq-d2q9-fluide} 
  \left \{ \begin {array}{l} \displaystyle 
f_0^{\rm eq}  =   {{\rho} \over{9}}    \, [ 1-\alpha +\beta \, ]  \,,  \\ \displaystyle  \vspace{-5 mm} \\ \displaystyle
f_1^{\rm eq}  =   {{\rho} \over{36}} \, [ 4 -\alpha-2\, \beta   
+ {{12  \, u_x }\over{\lambda}}   ] \,,\,\, 
f_2^{\rm eq}  =   {{\rho} \over{36}} \, \, \,[ 4 -\alpha-2 \,\beta   
+ {{12  \, u_y }\over{\lambda}}   ] \,, \\  \displaystyle  \vspace{-4 mm} \\ \displaystyle
f_3^{\rm eq}  =   {{\rho} \over{36}} \, [ 4 -\alpha-2 \, \beta   
- {{12  \, u_x }\over{\lambda}} ] \,,\,\, 
f_4^{\rm eq}  =   {{\rho} \over{36}} \, [ 4 -\alpha-2 \, \beta   
- {{12  \, u_y }\over{\lambda}} ] \,, \\  \displaystyle  \vspace{-4 mm} \\ \displaystyle
f_5^{\rm eq}  =   {{\rho} \over{36}} \, [ 4 +2 \, \alpha+\beta   
+{{3}\over{\lambda}}   (u_x+u_y  )  ] \,,\,\, 
f_6^{\rm eq}  =   {{\rho} \over{36}} \, [ 4 +2 \, \alpha+\beta   
+{{3}\over{\lambda}}   (-u_x+u_y  )   ] \,, \\  \displaystyle  \vspace{-4 mm} \\ \displaystyle
f_7^{\rm eq}  =   {{\rho} \over{36}} \, [ 4 +2 \, \alpha+\beta   
+{{3}\over{\lambda}}   (-u_x-u_y  ) ] \,,\,\, 
f_8^{\rm eq}  =   {{\rho} \over{36}} \, \, \,[ 4 +2 \, \alpha+\beta   
+{{3}\over{\lambda}}   (u_x-u_y  )   ] \, . 
 \end{array} \right.  
\monend 
In the space of moments, the linear collision is still described by a collision matrix.
In the fluid case, we have 
\moneq   \label{matrice-collision-fluide}
\!\!\! C \!=\! \left (
\begin{array}{ccccccccc}
\!\!\!\! 1 & \!\!\!\!\! 0 & \!\!\! 0 & \!\!\! 0 &  \!\!\! 0 &   \!\!\!\!\! 0 &  \!\!\!\!\! 0 &
\!\!\!\! 0\!\!\!\! &\!\!\!\! 0 \!\!\!\! \\ 
\!\!\!\!\!  0 & \!\!\!\!\! 1  
&  \!\!\! 0 &  \!\!\! 0 &  \!\!\! 0 &   \!\!\!\!\! 0 &  \!\!\!\!\! 0 &
\!\!\!\! 0\!\!\!\! &\!\!\!\! 0 \!\!\!\! \\ 
\!\!\!\!\!  0 & \!\!\!\!\! 0 &  \!\!\!  1  
 &  \!\!\! 0 &  \!\!\! 0 &   \!\!\!\!\! 0 &  \!\!\!\!\! 0 &
\!\!\!\! 0\!\!\!\! &\!\!\!\! 0 \!\!\!\! \\ 
\!\!\!\!\!  \alpha \, s_e \, \lambda^2 & \!\!\!\!\! 0 &  \!\!\! 0 &  \!\!\!\!\! 1-s_e &  \!\!\!\!\! 0 &  
 \!\!\! 0 &  \!\!\! 0 & \!\!\!\! 0 \!\!\!\!&\!\!\!\! 0 \!\!\!\! \\ 
\!\!\!\!\!  0 & \!\!\! 0 &  \!\!\!\! 0 &  \!\!\! 0 &  \!\!\! 1-s_x  &   \!\!\!\!\! 0 &  \!\!\!\!\! 0 
&\!\!\!\! 0 \!\!\!\!&\!\!\!\! 0 \!\!\!\! \\ 
\!\!\!\!\!  0 & \!\!\! 0 &  \!\!\!\! 0 &  \!\!\! 0 &  \!\!\! 0  &    \!\!\!\!\! 1-s_x &  \!\!\!\!\! 0 
&\!\!\!\! 0 \!\!\!\!&\!\!\!\! 0  \!\!\!\! \\ 
\!\!\!\!\!  0 & \!\!\!\!\! -s_q \, \lambda^2 
  &  \!\!\! 0 &  \!\!\! 0 &  \!\!\!\!\! 0 &   \!\!\!\!\! 0 & 
\!\!\!\!1-s_q &   \!\!\! 0 \!\!\!\!&\!\!\!\! 0 \!\!\!\! \\ 
\!\!\!\!\!  0 & \!\!\!\!\! 0 &  \!\!\!  -s_q \, \lambda^2  
 &   \!\!\! 0 &  \!\!\!\!\! 0 &   \!\!\!\!\! 0 &  
\!\!\! 0 & \!\!\!\!1-s_q \!\!\!\!&\!\!\!\! 0 \!\!\!\! \\ 
 \beta \, s_d \, \lambda^4 & \!\!\!\!\! 0 &  \!\!\! 0 &  \!\!\! 0 &  \!\!\! 0 &   
\!\!\!\!\! 0 &  \!\!\!\!\! 0 &\!\!\!\! 0\!\!\!\! &  \!\!\!\! 1-s_d 
\end{array} \right) \, . \monend
The sound velocity satisfies
\moneqstar  
c_0^2 = {{\alpha + 4}\over{6}} \, \lambda^2 \, . 
\monendstar 
The shear and bulk kinematic viscosities  $ \, \nu \, $ and $ \,  \zeta \, $ are given \cite{LL00} 
according to  the relations 
\moneqstar
\nu = {{\lambda}\over{3}} \, \Delta x  \, \Big( {{1}\over{s_x}} - {1\over2} \Big) \,, \quad  
\zeta =- {{\lambda}\alpha \over{6}} \, \Delta x  \, \Big( {{1}\over{s_e}} - {1\over2} \Big) \, . 
\monendstar
Note the usual values of the     
parameters~:   $ \, \alpha = -2 \,  $ and $ \, \beta = 1 $. 
During the relaxation step, the conserved variables 
$ \, W \equiv ( \rho \,,\, J_x  \,,\, J_y ) \,$ are    not   modified. 
The non-conserved moments  follow a relaxation algorithm 
described in (\ref{relaxation-moments}).

 \bigskip \bigskip   \noindent {\bf \large    3) \quad  
  Construction of the anti bounce back boundary condition }   

\begin{figure}    [H]  \centering 
  \centerline { \includegraphics[width=.50 \textwidth, angle=0] {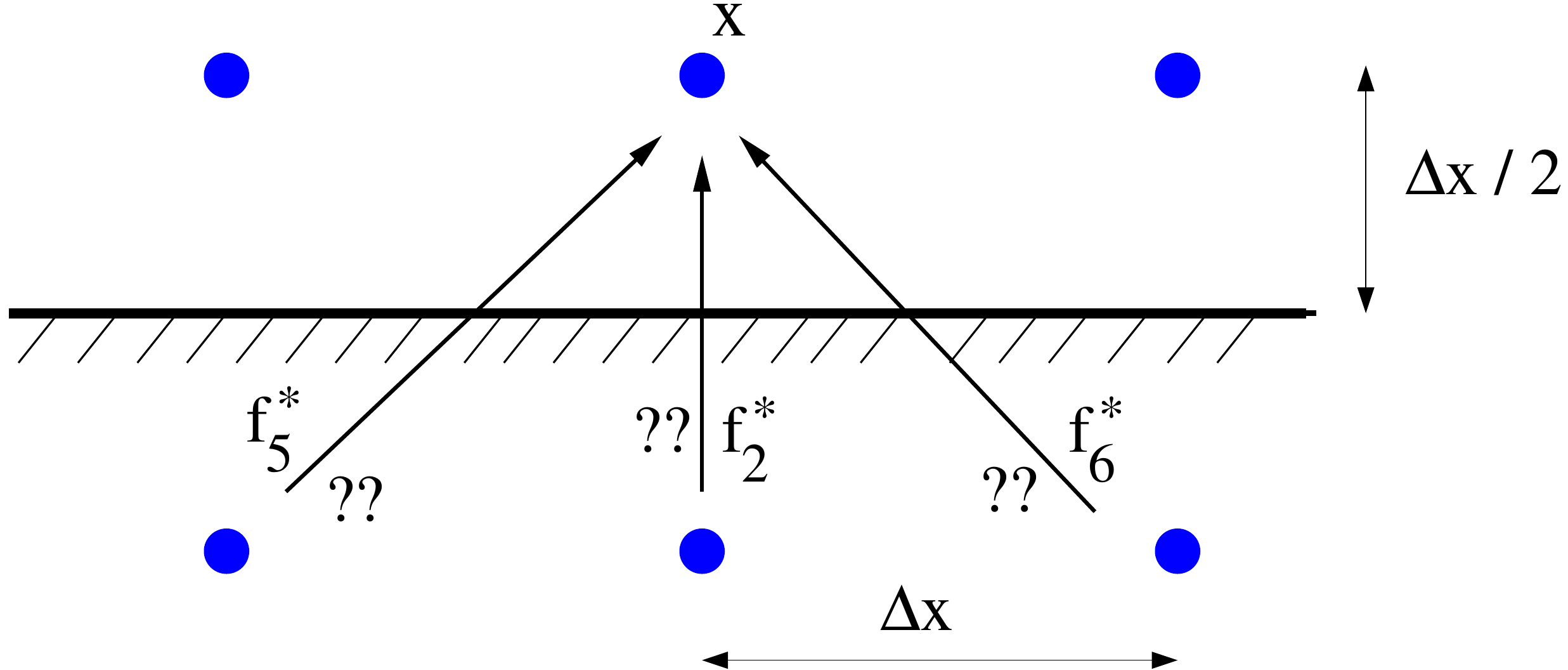}} 
\caption{The issue of choosing a boundary algorithm for the D2Q9 lattice Boltzmann scheme} 
\label{cl-d2q9-bas} \end{figure}

\noindent 
Consider to fix the ideas a bottom boundary for the D2Q9 lattice Boltzmann scheme, as illustrated in 
Figure \ref{cl-d2q9-bas}. For an internal node $ \, x $, 
the evolution of populations $ \, f_j \, $ is given by the internal scheme (\ref{lb-iteration}). 
For a bottom boundary, the  values of $ \, f_j^{*}(x-v_j \Delta t) \,$ 
for $ \, j  \in \{2,5,6\}\equiv \mathcal{B} \, $  are  {\it a priori} unknown.
%
%
When the vertex $ \, x \, $ is near the boundary, and located half a mesh size from 
the physical boundary,   
the internal iteration (\ref{lb-iteration}) is still applicable for the 
discrete velocities with numbers $ \, j = 0, \, 1, \, 3 ,\, 4, \, 7 ,\, 8 $, 
{\it id est} for $ \, j \not \in  \mathcal{B} $. 
For  $ \, j  \in  \mathcal{B} $, 
the fundamental scheme  (\ref{lb-iteration}) has to be adapted.
Our framework is still to try to  apply the internal scheme  at the  boundary.   
In our case, we set  
\moneqstar  \left \{ \begin{array} {l}
f_{5}(x, t+ \Delta t) =  f^*_{5}(x-(\Delta x,\, \Delta x), t ) \,, \\  
f_{2}(x, t+ \Delta t) =  f^*_{2}(x-(0,\, \Delta x), t )  \,, \\ 
f_{6}(x, t+ \Delta t)  =  f^*_{6}(x+(\Delta x,\, -\Delta x), t ) \,. 
\end{array} \right.  \monendstar 
Therefore, the boundary scheme is replaced by an extrapolation problem:
how to determine the values $ \,  f_j^{*}(x-v_j \, \Delta t,\, t) \,$ 
for  $ \, j  \in  \mathcal{B} \,$  from the known values  $ \,  f_k^{*}(x,\, t) \,$  
at the given vertex~$ \, x \, $ and the extra information given by 
the boundary conditions imposed by the physics?

\begin{figure}    [H]  \centering 
\centerline { \includegraphics[width=.35 \textwidth, angle=0] {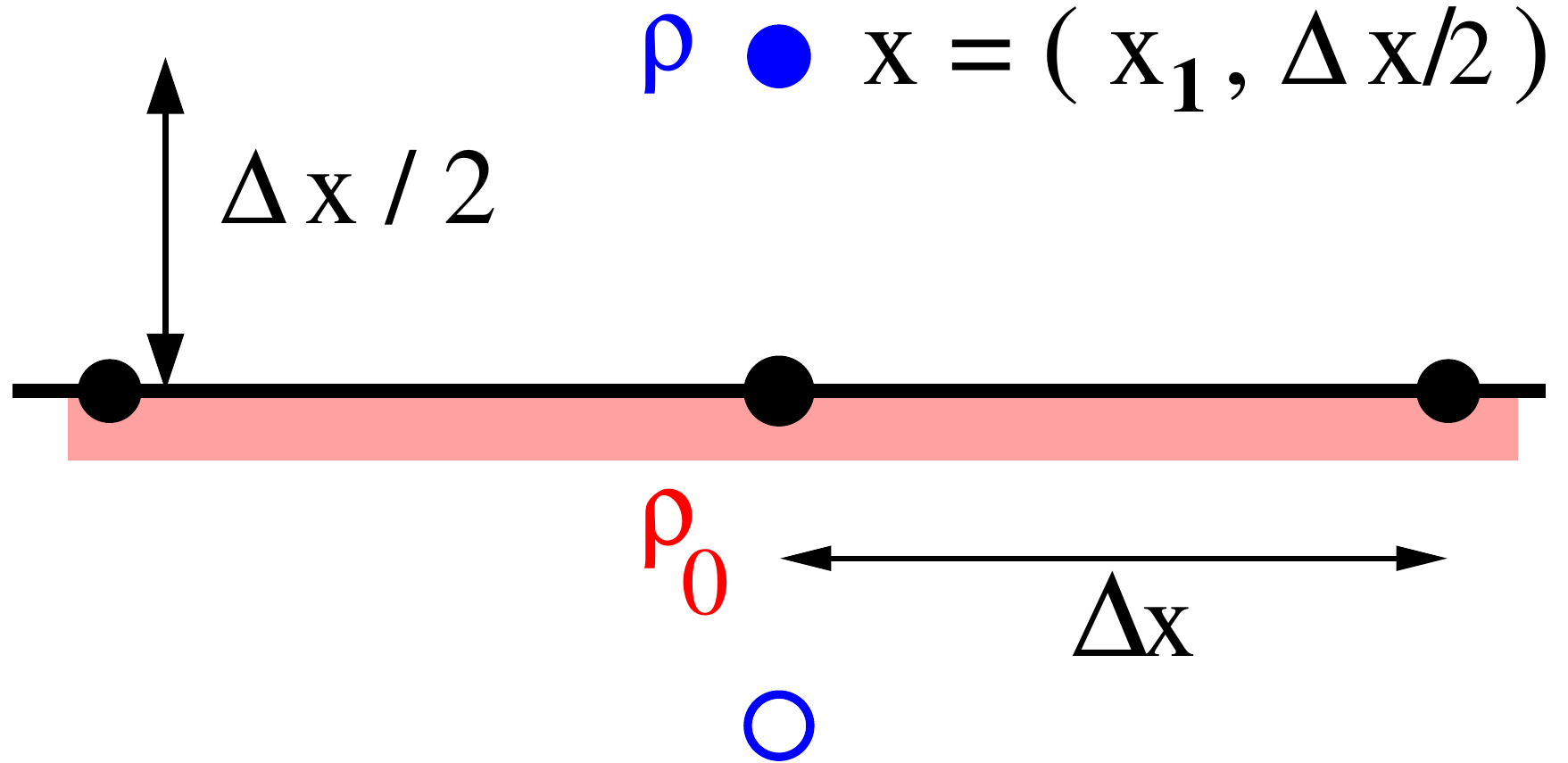}} 
\caption{Cell center framework for the boundary conditions: given  boundary value $ \, \rho_0 $ 
 and conserved moment $ \, \rho \, $ in the flow at vertex $ \, x  $.} 
\label{cl-d2q9-bas-data} \end{figure}

\monitem 
For the thermal case, the construction of the lattice Boltzmann scheme is based
on two arguments.  First we can approximate at first order of accuracy    
(see {\it e.g.} \cite{Du08}) 
the  values $ \,  f_j^{*}(x-v_j \Delta t,\, t) \,$ at the neighbouring vertices 
by the equilibrium values $ \,  f_j^{\rm eq}(x,\, t) \,$ at the internal vertex. 
Second, we have from the relations (\ref{feq-d2q9-thermique}) the following three   elementary remarks:
$ \, f_2^{\rm eq} =   {{\rho} \over{36}} \, ( 4 -\alpha-2 \,\beta ) $, 
 $ \, f_4^{\rm eq} =  {{\rho} \over{36}} \, ( 4 -\alpha-2 \, \beta  ) \, $ and in consequence 
$ \,  f_2^{\rm eq} +  f_4^{\rm eq} = {{4 -\alpha-2 \, \beta}\over{18}} \, \rho $.  
In a similar way, 
$ \, f_5^{\rm eq} = {{\rho} \over{36}} \, ( 4 +2 \, \alpha+\beta   ) \,  $ and 
$ \, f_7^{\rm eq} =  {{\rho} \over{36}}  \, (  4 +2 \, \alpha+\beta  ) \, $ implies 
$ \, f_5^{\rm eq} +  f_7^{\rm eq} = {{4 + 2 \, \alpha + \beta}\over{18}} \, \rho $. 
Last but not least, the relations 
$ \, f_6^{\rm eq} = {{\rho} \over{36}}  \, ( 4 +2 \, \alpha+\beta ) \, $
and 
$ \, f_8^{\rm eq} =   {{\rho} \over{36}} \, ( 4 +2 \, \alpha+\beta ) \, $
show  that 
$ \,  f_6^{\rm eq} +  f_8^{\rm eq} = {{4 + 2 \, \alpha + \beta}\over{18}} \, \rho $. 
If we suppose a Dirichlet boundary condition for the heat equation (\ref{chaleur}) then 
 the conserved ``temperature'' $ \, \rho = \rho_0 \, $ is  given on the boundary. 
The anti bounce back boundary condition  is obtained by replacing in the previous relations 
the equilibrium values by the outgoing particle distribution and the moment $ \, \rho (x) \, $    
by the given value $ \, \rho_0 \, $ on the boundary, as illustrated in Figure~\ref{cl-d2q9-bas-data}. 
%
\begin{figure}    [H]  \centering 
\centerline { \includegraphics[width=.50 \textwidth, angle=0] {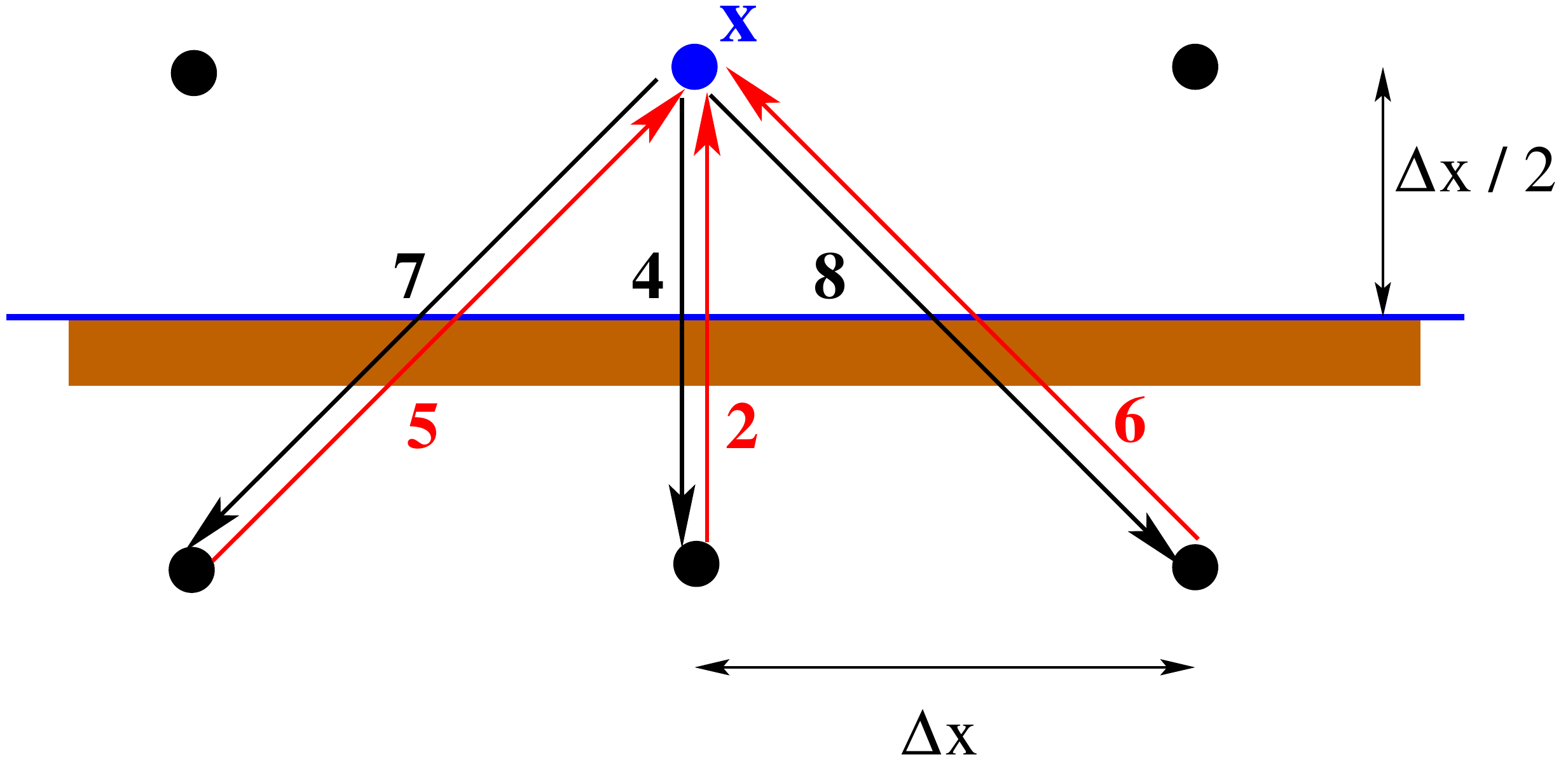}} 
\caption{For a boundary node $\, x \, $ along an horizontal frontier,
the D2Q9 lattice Boltzmann scheme   takes into account the three
numbers  $ \, \rho_0(x - {{\Delta x}\over{2}},\, t) $, 
$ \, \rho_0(x,\, t) \, $ and $ \, \rho_0(x + {{\Delta x}\over{2}},\, t) $.  } 
\label{cl-d2q9-oblique} \end{figure}
%
The  fluid node $ \, x \, $ near the boundary  has coordinates equal to 
$ \, (x_1,  {{\Delta x}\over{2}}) \,$  as depicted in Figure~\ref{cl-d2q9-bas-data}. 
In order to be more precise, we observe (see Figure \ref {cl-d2q9-oblique}) that 
the diagonal particles (numbered by $ \, j = 5 \, $ and $ \, j = 6 $) cross     
the boundary at locations $ \, x_1 - {{\Delta x}\over{2}} \, $ and 
$ \, x_1 + {{\Delta x}\over{2}} \, $ respectively. If we have a continuous information
for the given field on the boundary, we can introduce in the boundary scheme 
the exact values $ \, \rho_0(x_1 - {{\Delta x}\over{2}},\, t) $, 
 $ \, \rho_0(x_1 ,\, t) \, $ and $ \, \rho_0(x_1 + {{\Delta x}\over{2}},\, t) \, $
of the given data. 
Taking  into account all the previous remarks, the  anti bounce back boundary condition
is written in this case
\moneq \label{abb-chaleur}  \left \{ \begin{array} {l} \displaystyle 
f_{5}(x, \, t+ \Delta t) = -   f^*_{7}(x,\, t) \,+\,   {{4 + 2 \, \alpha + \beta}\over{18}} \, 
\rho_0 \Big(x_1 -\frac{\Delta x}{2}, \, t \Big)  \,, 
\\  \vspace{-4 mm} \\ \displaystyle
f_{2}(x, \, t+ \Delta t) =    -    f^*_{4}(x,\, t) \,+\, {{4 -\alpha-2 \, \beta}\over{18}} \, 
\rho_0 \big(x_1 , \, t \big)  \,,  
\\  \vspace{-4 mm} \\ \displaystyle 
\displaystyle f_{6}(x, \, t+ \Delta t)  =  -     f^*_{8}(x,\, t) \,+\,  {{4 + 2 \, \alpha + \beta}\over{18}} \, 
\rho_0 \Big(x_1 +\frac{\Delta x}{2}, \, t \Big) \, . 
\end{array} \right.  \monend
The minus sign in the right hand side of the relations (\ref{abb-chaleur}) characterizes the
denomination of ``anti'' bounce back.

\monitem 
For the linearized fluid, we have three conserved moments: density $ \, \rho \, $
and the two components of the momentum $ \, J \equiv \rho \, u $. 
The equilibrium of the particle distribution
is now given by the relation (\ref{feq-d2q9-fluide}). 
We have as previously three simple remarks.    
First, $\,  f_2^{\rm eq} = {{\rho} \over{36}} \, \big( 4 -\alpha-2 \,\beta   
+ {{12  \, u_y }\over{\lambda}} \big) \, $  
and $ \,  f_4^{\rm eq} =   {{\rho} \over{36}} \, \big( 4 -\alpha-2 \, \beta   
- {{12  \, u_y }\over{\lambda}} \big) $. In consequence, 
$ \,   f_2^{\rm eq} +  f_4^{\rm eq} = {{4 -\alpha-2 \, \beta}\over{18}} \, \rho $.  
Secondly, 
$ \, f_5^{\rm eq} =  {{\rho} \over{36}} \, \big( 4 +2 \, \alpha+\beta   
+{{3}\over{\lambda}}  (u_x+u_y ) \big) \, $ and 
$ \, f_7^{\rm eq} =   {{\rho} \over{36}} \, \big( 4 +2 \, \alpha+\beta   
+{{3}\over{\lambda}}  (-u_x-u_y ) \big) $.
After a simple addition of these two expressions, 
$ \,  f_5^{\rm eq} +  f_7^{\rm eq} = {{4 + 2 \, \alpha + \beta}\over{18}} \, \rho $.
Finally, the relations 
$ \, f_6^{\rm eq} =  {{\rho} \over{36}}  \, \big( 4 +2 \, \alpha+\beta   
+{{3}\over{\lambda}}  (-u_x+u_y )   \big) \,  $
and $ \,  f_8^{\rm eq} =  {{\rho} \over{36}} \, \big( 4 +2 \, \alpha+\beta   
+{{3}\over{\lambda}}  (u_x-u_y ) \big)   \,  $ 
show that 
$ \,   f_6^{\rm eq} +  f_8^{\rm eq} = {{4 + 2 \, \alpha + \beta}\over{18}} \, \rho $. 
%
Applying the same method as in the case of the thermal variant of the 
D2Q9 lattice Boltzmann scheme,  we can write the anti bounce back boundary condition
again with the relations  (\ref{abb-chaleur}). 
The interpretation of the variable $ \, \rho \, $ is now the density. 
The  anti bounce back is associated to a pressure datum, 
thanks to the acoustic hypothesis that $ \, p = c_0^2 \, \rho $. 
Observe that in the fluid case, the bounce back boundary condition is obtained simply by reversing the signs:
$ \, f_2^{\rm eq} -  f_4^{\rm eq} =  {{2}\over{3 \, \lambda}} \, \rho \, u_y  $, 
$ \, f_5^{\rm eq} -  f_7^{\rm eq} =  {{\rho}\over{6 \, \lambda}} (u_x+u_y ) \, $ and  
$ \, f_6^{\rm eq} -  f_8^{\rm eq} =  {{\rho}\over{6 \, \lambda}} (-u_x+u_y ) $. 

 \bigskip \bigskip   \noindent {\bf \large    4) \quad  
Linear asymptotic analysis 
for the heat problem }   

\begin{figure}    [H]  \centering 
  \centerline { \includegraphics[width=.35 \textwidth, angle=0] {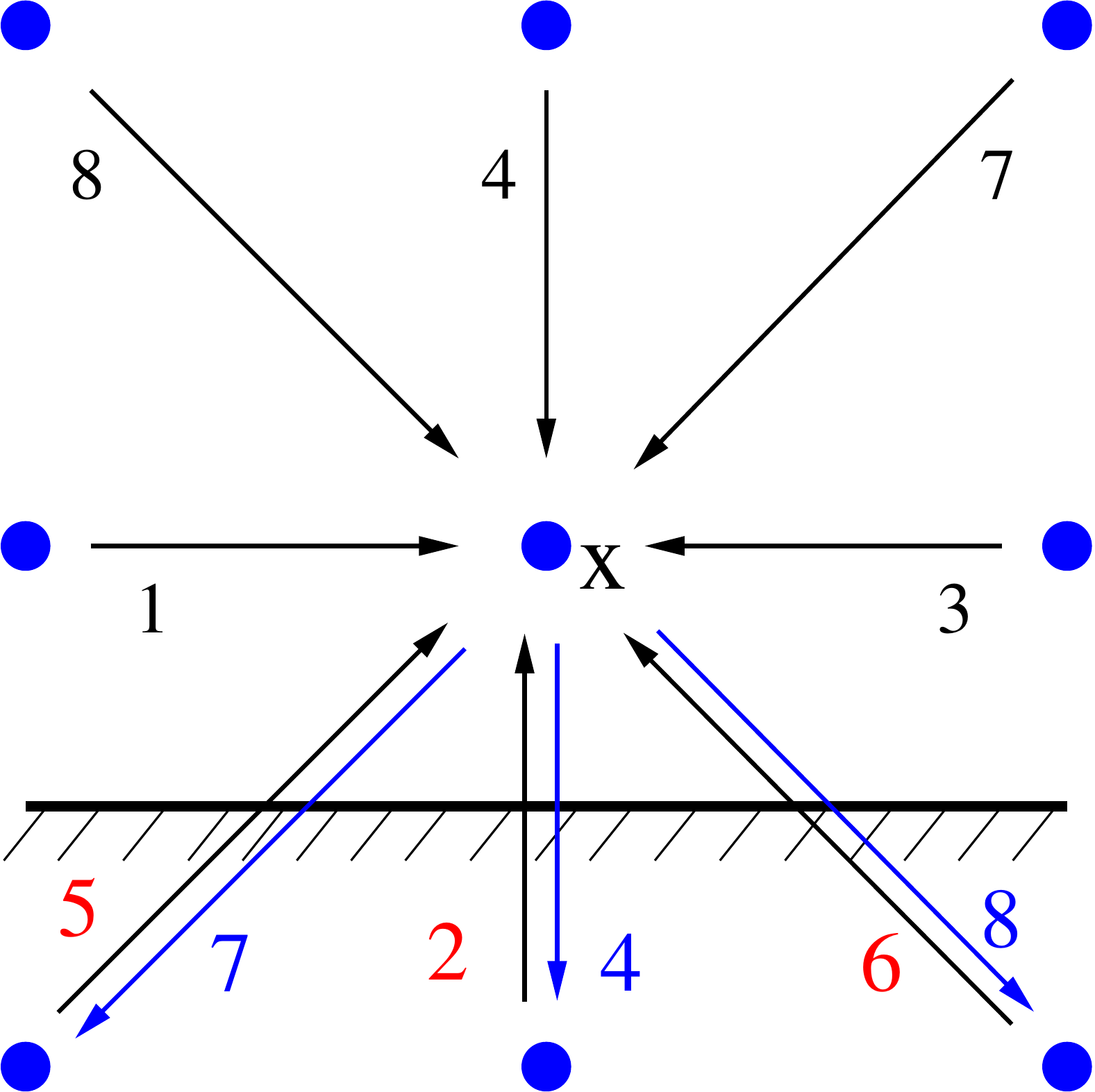}}
\caption{Vertex $ \, x \, $ near a boundary for the D2Q9 lattice Boltzmann scheme  } 
\label{cl-d2q9-global} \end{figure}

\noindent  
In order to analyze the  anti bounce back boundary condition  (\ref{abb-chaleur}), 
we can rewrite this condition as 
\moneqstar 
f^*_j(x, \, t ) =  - f^*_\ell(x,t)+\xi_j (x_j',t)  \, , \quad j \in \mathcal{B} \,,   
\monendstar 
with $ \,  \mathcal{B} = \{2, \, 5, \, 6 \} \, $ in our example. As previously 
the notation $ \, \ell \, $  corresponds to  the  opposite direction of the direction number $ \, j $:   $ \, v_j+v_\ell=0 $.  
%
The expression $ \, \xi_j(x_j',t) \, $  denotes the given  ``temperature''~$ \, \rho_0 \,$  on the boundary
at the space location $ \, x_j' $. 
In the specific example considered here, we have  
\moneq \label{cl-les-3-xi} \left \{ \begin{array} {l} \displaystyle 
 \xi_2  (x_2', \, t) =  {{4 -\alpha-2 \, \beta}\over{18}} \,  \rho_0 \big(x_1 , \, t)  \,, 
\\  \vspace{-4 mm} \\ \displaystyle 
 \xi_5 (x_5', \, t)  =   {{4 + 2 \, \alpha + \beta}\over{18}} \, \rho_0 \Big(x_1-\frac{\Delta x}{2}, \, t \Big) \,,\,
\\  \vspace{-4 mm} \\ \displaystyle 
 \xi_6 (x_6', \, t) 
=  {{4 + 2 \, \alpha + \beta}\over{18}} \, \rho_0  \Big(x_1+\frac{\Delta x}{2}, \, t \Big) \, . 
 \end{array} \right. \monend 
 If $ \, j \notin \mathcal{B} \, $ the above equation is replaced by the internal scheme: 
\moneqstar 
 f_j(x,t+\Delta t)=f^*_j(x - v_j\Delta t,t)  \,, \quad 
j = 0, \, 1, \, 3, \, 4 ,\, 7, \, 8  \, . 
\monendstar 
We have proposed in \cite{DLT15}   
a unified expression  of the lattice Boltzmann scheme D2Q9  for a  vertex   $ \, x \, $
near the boundary. With the help of  matrices  $ \, T_{j,\ell} \, $ and $ \, U_{j,\ell} $,
we can write
\moneq \label{cl-generale} 
f_j(x,t+\Delta t) =\sum_\ell T_{j,\ell} \, f_\ell^{*} (x,t) + \sum_\ell U_{j,\ell} \, f^*_\ell(x-v_j\Delta t,t) 
+ \xi_j (x_j', \, t) \, , \,\, 0 \leq j \leq 8 \, .  
\monend
The matrix element $ \, U_{j,\ell} \, $ is equal to $ \, 1 \, $ if $ \, \ell=j \notin \mathcal{B} \, $ 
and $ \, U_{j,\ell} = 0 \, $   if not:
\moneq \label{UU-cl-abb-chaleur} 
 U = \left( 
\begin{array}{ccccccccc}
{1} & 0 & 0 & 0 & 0 & 0 & 0 & 0 & 0 \\ 
0 & {1} & 0 & 0 & 0 & 0 & 0 & 0 & 0 \\ 
0 & 0 & 0 & 0 & 0 & 0 & 0 & 0 & 0 \\ 
0 & 0 & 0 & {1}  & 0 & 0 & 0 & 0 & 0 \\ 
0 & 0 & 0 & 0 & {1}  & 0 & 0 & 0 & 0 \\ 
0 & 0 & 0 & 0 & 0 & 0 & 0 & 0 & 0 \\ 
0 & 0 & 0 & 0 & 0 & 0 & 0 & 0 & 0 \\ 
0 & 0 & 0 & 0 & 0 & 0 & 0 & {1}  & 0 \\ 
0 & 0 & 0 & 0 & 0 & 0 & 0 & 0 &  {1}  
\end{array} \right) \, . \monend
In an analogous way, the explicitation of the opposite directions 
$ \, (5,\, 7) $,  $ \, (2, \, 4) \, $  and $ \, (6,\, 8 ) \, $   
leads to the following matrix denoted by $ \, T _, $ in the relation  (\ref{cl-generale}): 
\moneq \label{TT-cl-abb-chaleur} 
T = \left( 
\begin{array}{ccccccccc}
0 & 0 & 0 & 0 & 0 & 0 & 0 & 0 & 0 \\ 
0 & 0 & 0 & 0 & 0 & 0 & 0 & 0 & 0 \\ 
0 & 0 & 0 & 0 & \!  -  1 & 0 & 0 & 0 & 0 \\ 
0 & 0 & 0 & 0 & 0 & 0 & 0 & 0 & 0 \\ 
0 & 0 & 0 & 0 & 0 & 0 & 0 & 0 & 0 \\ 
0 & 0 & 0 & 0 & 0 & 0 & 0 & \!  - 1 & 0 \\ 
0 & 0 & 0 & 0 & 0 & 0 & 0 & 0 & \! -  1 \\ 
0 & 0 & 0 & 0 & 0 & 0 & 0 & 0 & 0  \\ 
0 & 0 & 0 & 0 & 0 & 0 & 0 & 0 & 0 
\end{array} \right)  \, . \monend
A complete picture of the boundary scheme is presented  in  Figure \ref{cl-d2q9-global}. 
%
%

\monitem 
The formal asymptotic analysis for infinitesimal space step $ \, \Delta x \, $
and infinitesimal time step $ \, \Delta t \, $ 
is  conducted as follows. Multiply the relations (\ref{cl-generale}) by the matrix  
$ \, M \, $ of (\ref{matrice-M}) in order to introduce the moments (\ref{moments}).
Thus a  complete expression of the D2Q9 lattice Boltzmann scheme  for a node   $ \, x \, $    near the boundary    
takes the form
\moneqstar 
m_k(x,t+\Delta t) =  (M T M^{-1} C)_{k,\ell} \ m_\ell (x,t)  + 
 (M_{k,\ell} \, U_{\ell,j} \,  M^{-1}_{j,p} \,  C_{p,q}) \ m_q (x-v_\ell \Delta t,t) + 
 M_{k,\ell} \, \xi_\ell   \, . 
\monendstar  
Then  expand this relation  at order 0 and 1.  
At order  zero, we have 
\moneqstar 
m_k + {\rm{O}}(\Delta t)  = (M T M^{-1} C)_{k,\ell} \ m_\ell 
+ (M U M^{-1} C)_{k,\ell} \ m_\ell  + {\rm{O}}(\Delta x)  +  M_{k,\ell} \, \xi_\ell   \, . 
\monendstar   
Drop away the indices: 
\moneqstar 
m + {\rm{O}}(\Delta t)  = (M T M^{-1} C) \, \ m + (M U M^{-1} C) \, m  +  M \, \xi + {\rm{O}}(\Delta x)   \, . 
\monendstar   
We make a hypothesis of acoustic scaling: the ratio $\, {{\Delta x}\over{\Delta t}} \, $
is maintained constant as $ \, \Delta x \, $ and $ \, \Delta t \, $ tend to zero. 
We put also all the  unknowns in the left hand side:  
\moneqstar 
 \big( I - M (T+U) M^{-1} C \big) \, m =  M \, \xi + {\rm{O}}(\Delta x)  \, .  
\monendstar   
We introduce the matrix $ \, K \, $ according to 
\moneq \label{matrice-K}  
 K \equiv I - M (T+U) M^{-1} C  \, .
\monend
We have established  asymptotic equations for the moments $ \, m \, $ in the first cell
of  the computational domain for anti bounce back scheme  at  order zero:  
\moneq \label{abb-systeme-lineaire}  
  K \, m \,= \, M \, \xi + {\rm{O}}(\Delta x) \, . 
\monend
with
\moneqstar
\xi =  \xi_0   + \Delta t \, \partial \xi + {\rm{O}}(\Delta x^2) \, .  
\monendstar 

\monitem 
For the  anti bounce back for the thermal D$2$Q$9$ scheme, we have:
\moneq \label{abb-thermic-matrice-K}  
 { K=\left( \begin{array}{ccccccccc}
\!\!\!\! {1\over6} (4+\alpha s_e)  & \!\!\! 0 &  \!\!\!\!  0  &  \!\!\!\! {1\over{6 \lambda^2}} (1-s_e)  
&  \!\!\!\! {{1}\over{2 \lambda^2}} (s_x-1) &  \!\!\!\!\!\!\! 0 &  \!\!\!\! 0 &  \!\!\!\! 0 &  \!\!\!\! 0 \\ 
\!\!\!\! 0 &  \!\!\!\!  s_j  &  \!\!\!\!  0 &  \!\!\!\! 0 &  \!\!\!\! 0 & 
 \!\!\!\! {{1}\over{\lambda}} (1-s_x)   &   \!\!\!\! 0 &  \!\!\!\! 0 &  \!\!\!\! 0  \\ 
\!\!\!\!   {{\lambda}\over{6}} (4+\alpha s_e)  & \!\!\!\!  0 &  \!\!\!\!   s_j & 
 \!\!\!\!  {{1}\over{6  \lambda}} (1 - s_e) &  \!\!\!\!  {{1}\over{2 \, \lambda}}  (s_x - 1)
&  \!\!\!\!\!\!\! 0 &  \!\!\!\! 0 &   \!\!\!\! 0 &  \!\!\!\! 0  \\ 
\!\!\!\!  {{\lambda^2}\over{6}} (4 - 3 \alpha s_e + 2 \beta s_d )
 & \!\!\!\!  0 &  \!\!\!\!  0  &   \!\!\!\! {{1}\over{2}} (s_e + 1)
 &  \!\!\!\!    {{1}\over{2}} (1 - s_x) &  \!\!\!\!\!\!\! 0 &    \!\!\!\! 0 &  \!\!\!\! 0 & 
 \!\!\!\! {{1}\over{3 \lambda^2}} (1 - s_d) \\ \vspace {-.5 cm}  \\ 
\!\!\!\!   {{\lambda^2}\over{18}} (-4 + \alpha s_e + 2 \beta s_d )
 &  \!\!\!\!   0  &  \!\!\!\! 0 &  \!\!\!\!  {{1}\over{18}} (1 - s_e) &  \!\!\!\!  {{1}\over{2}} (1 + s_x)
 &  \!\!\!\!\!\!\! 0  &  \!\!\!\! 0  &  \!\!\!\! 0 &  \!\!\!\! {{1}\over{9 \lambda^2}}  (1 - s_d) \\ 
\!\!\!\! 0 & \!\!\!\!  0 &  \!\!\!\! 0 & \!\!\!\!  0 &  \!\!\!\! 0 &  \!\!\!\!\!\!\!1 &
  \!\!\!\!0 &  \!\!\!\!0 &  \!\!\!\!0 \\ 
\!\!\!\! 0 & \!\!\!\! 0  
&  \!\!\!\! 0 &  \!\!\!\! 0 &  \!\!\!\! 0 & \!\!\!\!\!\!\!\! \lambda (1 - s_x) & \!\!\!\! s_q & \!\!\!\! 0 & 0 \\ 
\!\!\!\!   {{\lambda^3}\over{3}}  (\alpha s_e + \beta s_d ) & \!\!\!\!  0 &  \!\!\!\!  0 &
 \!\!\!\!  {{\lambda}\over{3}} (1-s_e ) &  \!\!\!\! \lambda (1-s_x) &  \!\!\!\!\!\!\! 0 &  \!\!\!\! 0 &  \!\!\!\! s_q &
 \!\!\!\!  {{1}\over{3 \lambda}} (1-s_d) \\ \vspace {-.5 cm}  \\ 
\!\!\!\!   {{\lambda^4}\over{3}}  (\alpha s_e - 2 \beta s_d ) &  \!\!\!\!  0 &  \!\!\!\!  0 &
 \!\!\!\!  {{\lambda^2}\over{3}} (1-s_e) &   \!\!\!\! \lambda^2 (1-s_x) &  \!\!\!\!\!\!\! 0 &  \!\!\!\! 0 & 
 \!\!\!\! 0 &  \!\!\!\! {{1}\over{3}} (1 + 2 s_d)
\end{array} \right)} \, . \monend
The matrix $ \, K \, $ defined in (\ref{abb-thermic-matrice-K}) is   regular. 
The solution  of anti bounce back  at order zero $ \,  m_0  \, $ is the unique solution of the 
equation 
\moneq \label{abb-systeme-lineaire-zero}  
  K \, m_0 =  M \, \xi_0  \, 
\monend
obtained by neglecting the first order terms in (\ref{abb-systeme-lineaire}),
with $ \, \xi_0 \, $ given by the relations (\ref{cl-les-3-xi}):
%
%
%
\moneqstar 
\xi_0  =   \, \big( 0 \,,\, 0 \,,\, {1\over18} \,(-\alpha - 2\, \beta + 4) \, \rho_0 \,,\, 0 \,,\, 0 \,,\,
{1\over18} \,(2\alpha + \beta + 4) \, \rho_0 \,,\,{1\over18} \,(2\alpha + \beta + 4) \, \rho_0 \,,\,  0 \,,\, 0 \big)^{\textrm t} 
\monendstar 
Then  $ \,  m_0 = K^{-1} \,\smb \,  M \, \xi_0     \, $ 
since the matrix $ \, K \, $ is regular. 
 With the given temperature  $ \, \rho_0 \, $  given on the boundary, we have:
\moneq \label{abb-thermique-ordre-zero}  
  m_0=\left(\rho_0, \, 0, \,0, \,\alpha \, \rho_0  \, \lambda^2,\, 0,\, 0,\, 0,\, 0, \,
\beta \, \rho_0 \,  \lambda^4 \right )^{\rm t}  \, .  
\monend

\monitem  
At order one, the asymptotic analysis of anti bounce back  follows the work 
done in \cite{DLT15}  for the usual bounce back boundary condition. 
We introduce the matrices
%
\moneqstar 
 B_{k,p}^{\alpha} = \sum_{\ell,j,q} M_{k,\ell} \, U_{\ell,j}  \,v_j^{\alpha}  \,M_{j,q}^{-1}  \,C_{q,p}  \,, \,\,\,
\alpha = 1, \, 2 \, .  
\monendstar 
The equivalent equations for anti bounce back scheme  up to order one are solutions  of  
\moneqstar 
 K \, m=M \, \xi_0 + \Delta t \left[ M \, \partial \xi - 
\partial_t m - B^{\alpha} \, \partial_\alpha m \right] + {\rm{O}}(\Delta x^2)  \, .  
\monendstar 
Recall that with the data (\ref{cl-les-3-xi}), we have 
%
%
\moneqstar 
\partial \xi = \Big( 0 \,,\, 0\,,\, 0 \,,\, 0\,,\,  0 \,,\, {{\lambda}\over{36}} \, (-2\, \alpha - \beta - 4) \, \partial_x \rho_0
\,,\, {{\lambda}\over{36}} \, (2\, \alpha + \beta + 4) \, \partial_x \rho_0 \,,\,  0 \,,\, 0 \Big)^{\textrm t}  \, . 
\monendstar 
We search a formal  expansion of the type 
\moneqstar 
m=m_0  + \Delta t \,\,  m_1 + {\rm{O}}(\Delta x^2) \,. 
\monendstar 
We have  consequently
\moneqstar  \, m_1 = K^{-1} \, \smb \, \left( M  \, \partial \xi_0 
- \partial_t m_0 - B^{\alpha} \partial_\alpha m_0 \right) \, . 
\monendstar 
After a tedious computation, the conserved variable $\, \rho \, $  can finally be expanded as 
\moneq \label{abb-thermique-ordre-un}  
 \rho \Big(x_1, {{\Delta x}\over{2}} \Big) = \rho_0 (x_1)  + {1\over2} \, \Delta x \,\, \partial_y  \rho (x_1,\, 0) 
+  {\rm{O}}(\Delta x^2)  \, . 
\monend
The derivative $ \, \partial_y  \rho (x_1,\, 0) \, $ is the exact value for the continuous problem
at the boundary point. 
Finally, the relation (\ref{abb-thermique-ordre-un}) is nothing else than the Taylor expansion for 
the conserved variable between the boundary value  $ \, \rho_0 (x)\, $     
and the value $ \,  \rho (x) \, $ in the first fluid  vertex. 
This Taylor expansion is correctly recovered at first order.
This indicates that the Dirichlet boundary condition
is correctly taken into consideration and does not induce spurious effects  \cite{WWLL13}. 
%
Recall that for general bounce back velocity
condition, the situation shows several discrepancies  \cite{DLT15, DLT17}.  


\monitem  
In order to illustrate the good quality of the results, 
we have evaluated numerically the eigenmodes of the heat  equation   
$ \,\, \partial_t\rho-\mu \, \Delta\rho = 0 $, {\it id est}  the solutions of the
modal  problem  
\moneqstar 
-\Delta \rho = \Gamma \, \rho \quad
\monendstar 
in a square domain $ \, \Omega = \, ] 0,\, L[^2 \,$ with  
periodic and Dirichlet boundary conditions. The exact eigenvalues 
$ \, \Gamma_{k,\, \ell} \, $ are proportional to 
$ \, k^2 + \ell^2 $ (the index of the mode).   
For periodic  boundary condition, $\, k \, $ and $ \, \ell \, $ are 
even integers, then $ \,  k^2 + \ell^2 = 4,\, 8,\, 16,\, ...$ 
For homogeneous  Dirichlet boundary condition, there is no
restriction on the integers  $\, k \, $ and $ \, \ell $; then  
$ \,  k^2 + \ell^2 = 2,\, 5,\, 8,\, ...$ 
We use the D2Q9 lattice Boltzmann scheme with homogeneous anti
bounce back boundary condition on a $ \, 71 \times 71 \, $ grid. 
The boundary is always located exactly between two grid points 
and the previous scheme is used without any modification. 
The eigenmodes of the problem are evaluated in the following way~:   
after one step of the algorithm, the field $ \, f(x, \, t+\Delta t) \, $ is 
supposed to be proportional to  $ \, f(x, \, t) $. The modes are computed 
numerically with the  ``ARPACK'' Arnoldi package  \cite{arpack}. 

\begin{figure}    [H]  \centering 
\centerline { \includegraphics[width=.88 \textwidth, angle=0] {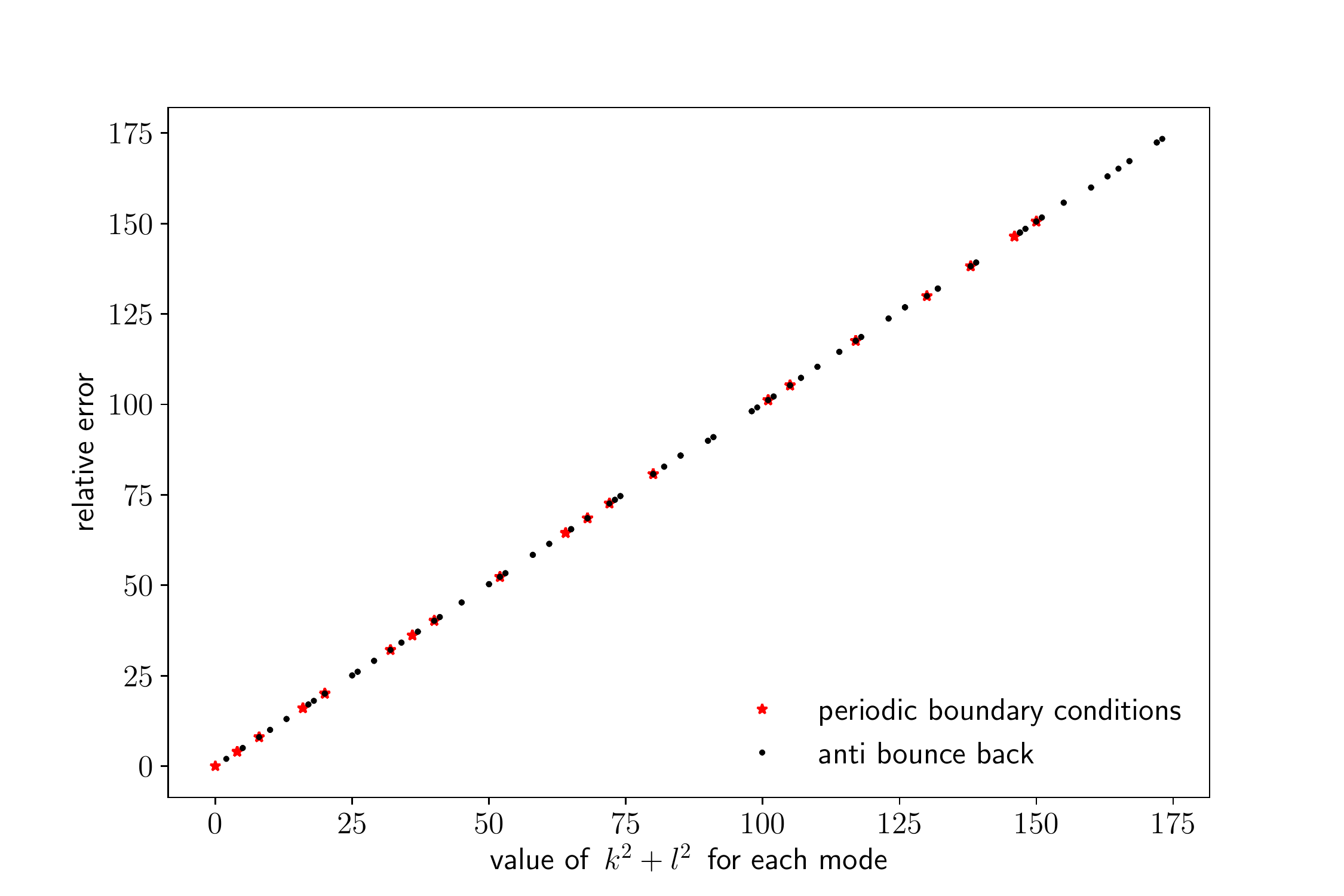}}
\caption{Relative error $ \, {{\Gamma_{\rm computed}}\over{\Gamma_{\rm exact}}} - 1 \, $ 
(multiplied by $10^4$) 
as a function of the index $ \,  k^2 + \ell^2 \, $ of the modes  
for the Laplace equation in a square. }    
\label{erreur-modes-laplace-carre} \end{figure}

\newpage \noindent 
 On Figure~\ref{erreur-modes-laplace-carre}, we have plotted the relative error (multiplied by $10^4$)
\moneqstar  
\varepsilon_{k,\, \ell,\, {\rm computed}} =  {{\Gamma_{k,\, \ell,\, {\rm computed}} \over{\Gamma_{k,\, \ell,\, {\rm exact}}}}} - 1 \, . 
\monendstar 
First, this error is very small: less than $ \, 0.025 \, $ \%. 
Second, it is approximately   
proportional to the  index $ \,  k^2 + \ell^2 $.
This indicates a second order accuracy with respect to the wave number     
$ \,  \pi \, {{\sqrt{ k^2 + \ell^2}}\over{L}} $, with $ \, L \, $ the size of the computational domain.
This second order error can be minimized by a suitable
choice (see {\it e.g.} our contribution  \cite{DL11})  of the parameters of the scheme
\moneqstar   
 s_j = 1.20 ,\,  s_e = 1.30 ,   \,  
  s_x = {{1}\over{{{\sqrt{3}}\over{6}} + {1\over2}}} \simeq 1.27 ,   \,  
 s_q = {{1}\over{{{\sqrt{3}}\over{3}} + {1\over2}}} \simeq 0.93 ,   \,  
  s_d = 1.70    ,   \,   \alpha = -2  ,   \,  \beta = 1 \, . 
\monendstar 
%
%
%

 \bigskip \bigskip   \noindent {\bf \large    5) \quad  
  Extended anti bounce back for general boundaries  }   

\begin{figure}    [H]
\vskip -1.5 cm
\centerline  { \quad \includegraphics[width=1.15 \textwidth, angle=0] {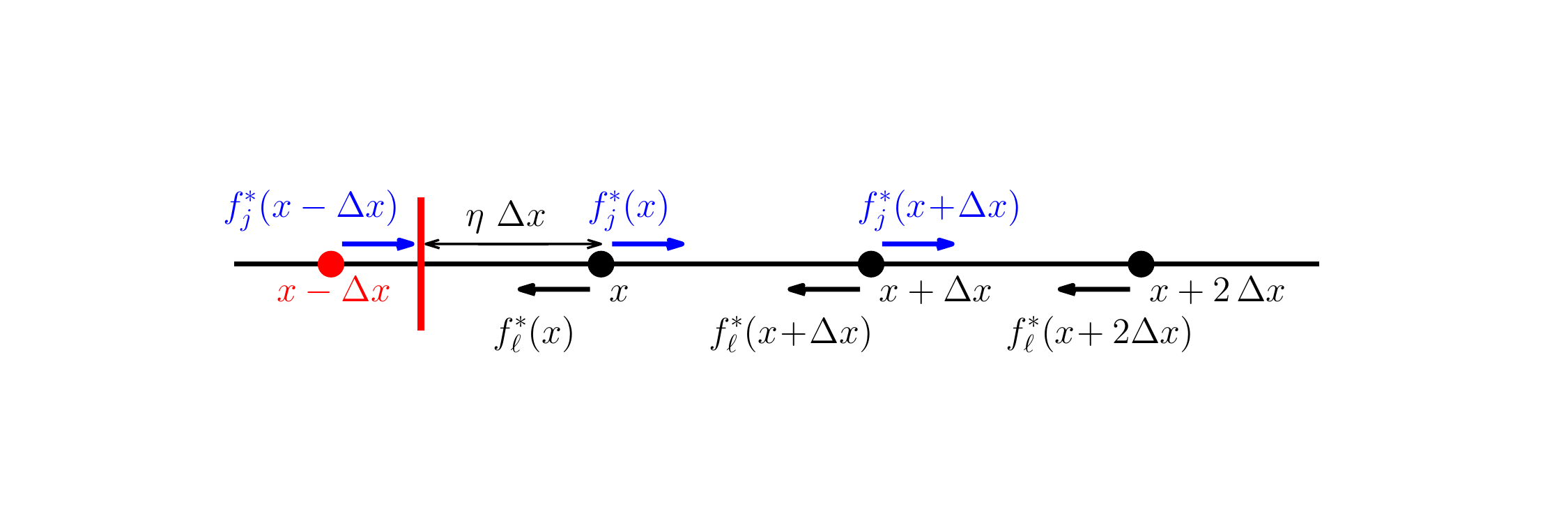}}
\vskip -1.5 cm
\caption{Boundary conditions when the boundary is not in the middle of two
mesh vertices} 
\label{cl-bfl-jpg} \end{figure}

\noindent  
The  argument summarized in the previous section for  bounce back  
and anti bounce back  conditions can be extended to consider a  
boundary at an arbitrary distance $ \, \eta  \, \Delta x \, $ ($ 0 < \eta < 1$) from
the last fluid point.
It is completely described for bounce back in Bouzidi {\it et al.} \cite{BFL01}. 
Assume that the boundary is on the left for a given space direction $ \, e_j \, $
directed towards in interior of the computational domain. 
 The boundary node $ \, x \, $ is located inside the fluid 
and the node $ \, x+\Delta x \, $ on its right is also inside the fluid.
But the node $ \, x-\Delta x  \, $ is not a fluid node.   
The outgoing particles $ \, f^*_{j}(x) \, $ and $ \, f^*_{j}(x+\Delta x ) \, $ 
are known at time $ \, t \, $ after the relaxation step. 
The incoming particle  $ \, f^*_\ell (x-\Delta x ) \, $ 
for the  opposite direction $ \, - v_j \, $ 
($  v_j+v_\ell=0 $)   has to be determined by the boundary condition. 
This question is summarized in Figure~\ref{cl-bfl-jpg}.    
%
%
Assuming a linear variation of the macroscopic properties of    
the flow, one gets (see the original derivation in \cite{BFL01})
extended  bounce back  to impose the velocity in the fluid case: 
%
%
\moneq \label{fbl01-bb}  
f^*_j (x-\Delta x) =  \left \{ \begin{array} {l} \displaystyle 
2 \, \eta\, f^*_\ell(x)+(1-2\, \eta) \, f^*_\ell(x+\Delta x)+ \xi_j \qquad  {\rm for} \, \eta \leq {1\over2}  
\\ \displaystyle 
\frac{1}{2\, \eta} \,  f^*_\ell(x) +  \Big( 1-\frac{1}{2\, \eta}  \Big) \, f^*_j(x)+ 
{{\xi_j}\over{2\, \eta}} \, \qquad \qquad {\rm for} \, \eta \geq {1\over2}     
\end{array} \right.  \monend
or extended anti bounce back (to impose the density in the fluid case):
%
%
\moneq \label{fbl01-abb}  
f^*_j (x-\Delta x) =  \left \{ \begin{array} {l} \displaystyle 
-2 \, \eta\, f^*_\ell(x) - (1-2\, \eta) \, f^*_\ell(x+\Delta x)+ \xi_j \qquad  {\rm for} \, \eta \leq {1\over2}  
\\ \displaystyle 
-\frac{1}{2\, \eta} \,  f^*_\ell(x) +  \Big( 1-\frac{1}{2\, \eta}  \Big) \, f^*_j(x)+ 
{{\xi_j}\over{2\, \eta}} \, \qquad \qquad {\rm for} \, \eta \geq {1\over2}  \, . 
\end{array} \right.  \monend
The interpolation is done before the propagation if $ \, \eta \leq {1\over2} \, $
and after if $ \, \eta \geq {1\over2} $.

\monitem
Similar expressions have been written  in \cite{BFL01} for extended  bounce back  assuming a parabolic variation of the
velocity field 
\moneq \label{fbl02-bb}  
  f^*_j (x-\Delta x) =   \left \{ \begin{array} {l} \displaystyle 
 \eta\, (1+2\, \eta) \, f^*_\ell(x) + (1-4\, \eta^2) \, f^*_\ell(x+\Delta x)
    \\  \displaystyle  \qquad \qquad \qquad   -  \eta \, (1 - 2 \, \eta) \, f^*_\ell(x+2\, \Delta x)
  + \xi_j \qquad  {\rm for} \, \eta \leq {1\over2}  
\\ \displaystyle \frac{1}{\eta \, (1+2\, \eta)} \, \Big(  f^*_\ell(x) + \xi_j  \Big) 
  +  \Big( 2 - {{1}\over{\eta}} \Big) \, f^*_j(x)
\\  \displaystyle \vspace {-4 mm}  \\  \displaystyle\qquad \qquad \qquad +  {{1 - 2 \, \eta}\over{1 + 2 \, \eta}} \,  f^*_j(x + \Delta x) 
 \, \qquad \qquad \qquad {\rm for} \, \eta \geq {1\over2}    \, .  
\end{array} \right.  \monend
  Similar changes in the signs as above define the extended for anti bounce back 
\moneq \label{fbl02-abb}  
  f^*_j (x-\Delta x) =   \left \{ \begin{array} {l} \displaystyle 
 -\eta\, (1+2\, \eta) \, f^*_\ell(x) - (1-4\, \eta^2) \, f^*_\ell(x+\Delta x)
    \\  \displaystyle  \qquad \qquad \qquad   +  \eta \, (1 - 2 \, \eta) \, f^*_\ell(x+2\, \Delta x)
  + \xi_j \qquad  {\rm for} \, \eta \leq {1\over2}  
\\ \displaystyle \frac{1}{\eta \, (1+2\, \eta)} \, \Big(  -f^*_\ell(x) + \xi_j  \Big) 
  +  \Big( 2 - {{1}\over{\eta}} \Big) \, f^*_j(x)
\\  \displaystyle \vspace {-4 mm}  \\  \displaystyle\qquad \qquad \qquad +  {{1 - 2 \, \eta}\over{1 + 2 \, \eta}} \,  f^*_j(x + \Delta x) 
 \, \qquad \qquad \qquad {\rm for} \, \eta \geq {1\over2}   \, .  
\end{array} \right.  \monend

\newpage 

\monitem    Application  with the heat equation 

\noindent 
The above anti bounce back boundary condition has been used in our previous contribution  \cite {DL09}
in the context of the heat equation.
Numerical experiments are presented in this contribution. Figures  3 and 4  with D2Q5 scheme for thermics inside a circle, 
and Figures 6 to 9 for D3Q7 lattice Boltzmann scheme for thermics in a  three-dimensional ball.
The Dirichlet boundary condition
is simply implemented with an extended anti bounce back scheme.  

\monitem    Applications with  linear acoustics.

\noindent 
For linear acoustics, we have used in  \cite {DL11} the time dependent  extended anti bounce back (\ref{fbl02-abb})
to enforce a given harmonic pressure on the boundary of a disc (see the Figures~1 to 3 of this reference). 
No particular difficulty was reported with this treatment of the
time dependent pressure boundary condition. 

\begin{figure}    [H]  \centering 
\vskip -1.8  cm
\centerline { \includegraphics[width=.95  \textwidth, angle=0] {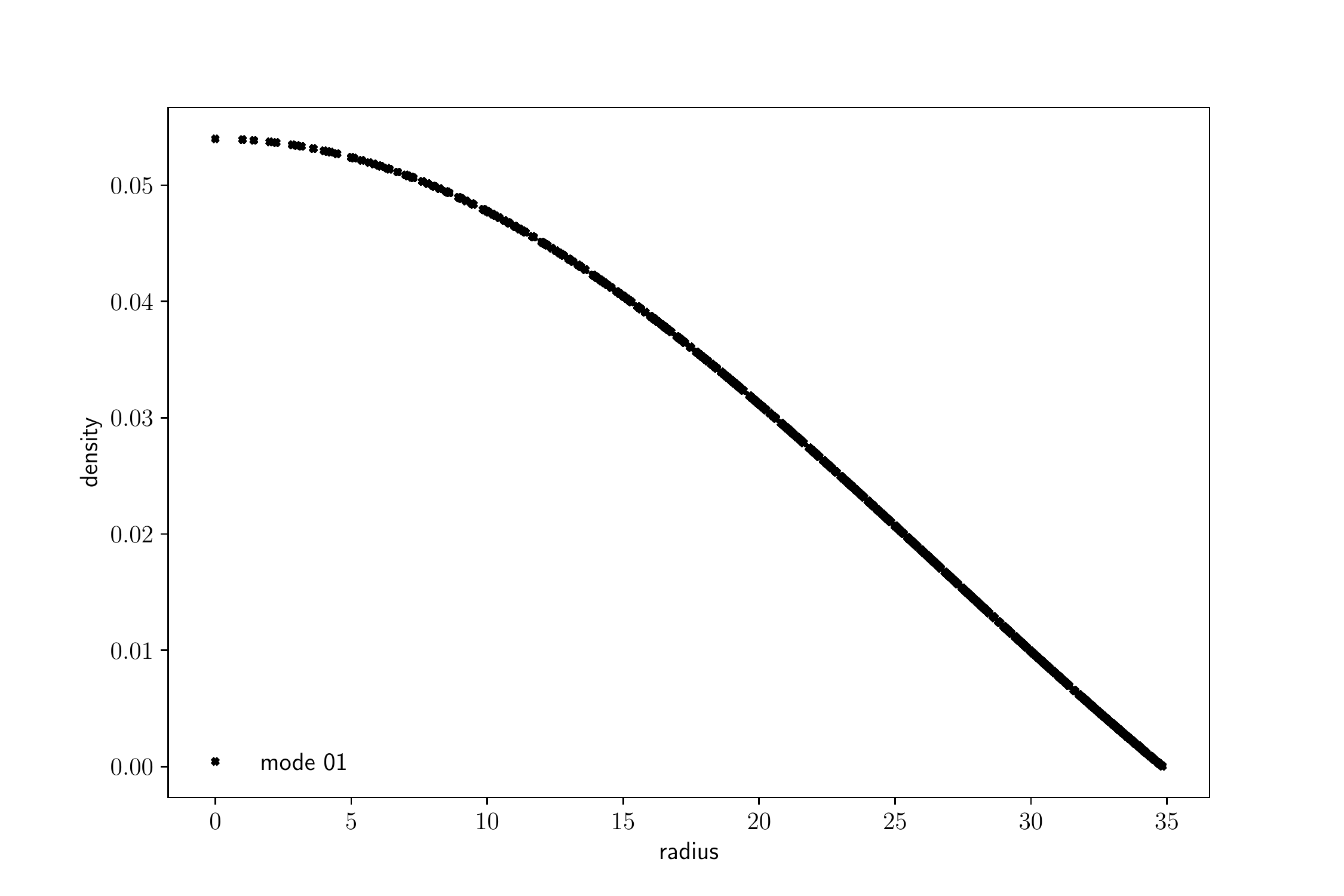}}
\vskip -.5  cm
\caption{Linear fluid problem in a disc. Density profile as a function of the radius for the first eigenmode
($\Gamma= -0.03983758 + 0.00019317 \,\, i \,$ in unit $ \,  1/\Delta t$) 
with homogeneous  anti bounce back boundary condition for a disc with a radius
$ \, R = 34.85 $. } 
\label{abb-mode-34-densite} \end{figure}

\begin{figure}    [H]  \centering 
\vskip -1.8  cm
\centerline { \includegraphics[width=.95  \textwidth, angle=0] {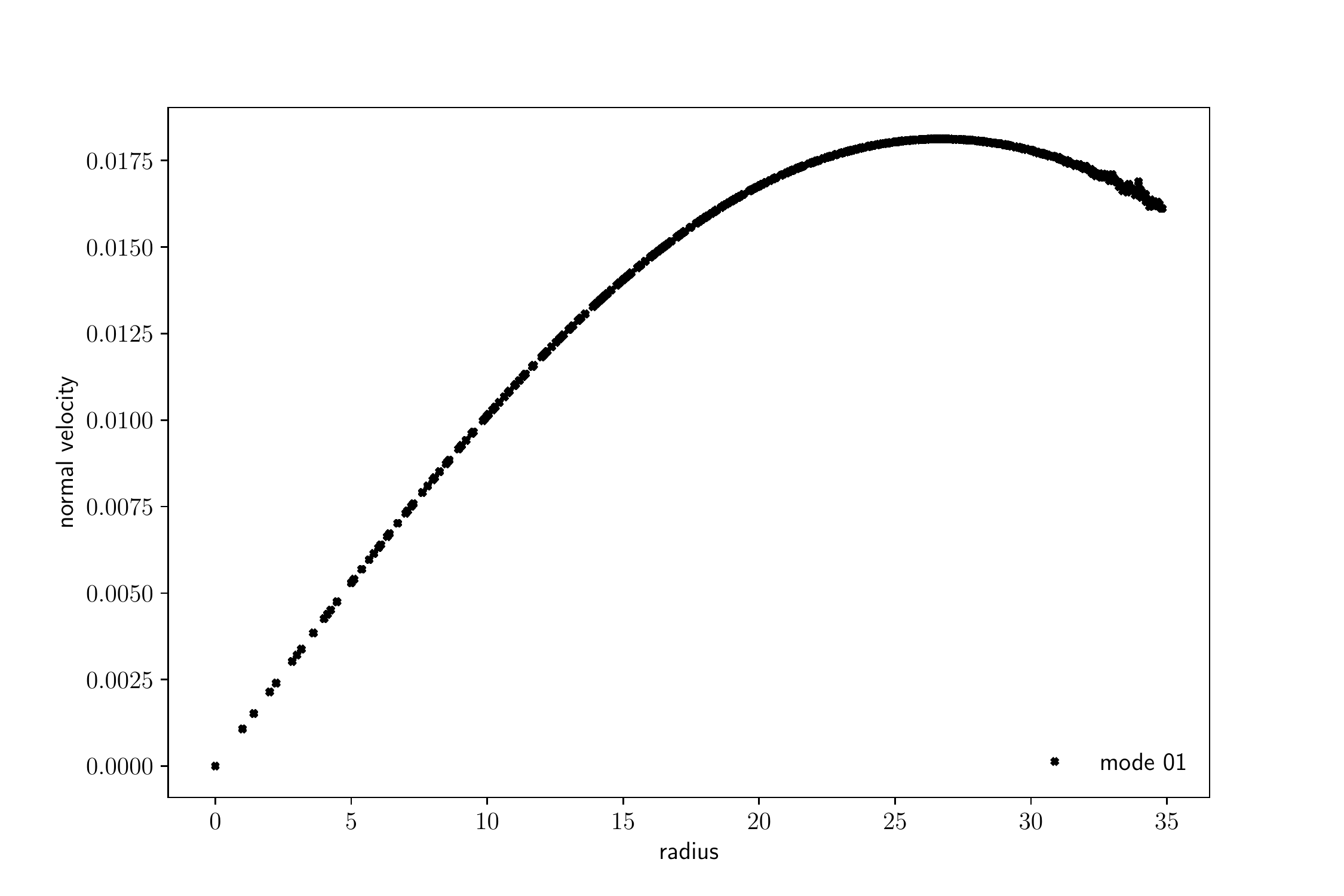}}
\vskip -.5  cm
\caption{Linear fluid problem in a disc. Radial velocity profile  as a function of the radius for the first eigenmode
($\Gamma= -0.03983758 + 0.00019317 \,\, i \,$ in unit $ \,  1/\Delta t$) 
with homogeneous  anti bounce back boundary condition for a disc with a radius
$ \, R = 34.85 $. A small discrepancy due to treatment of the
boundary condition is visible.   The tangential velocity is negligible.   } 
\label{abb-mode-34-vit-radiale} \end{figure}

\begin{figure}    [H]  \centering 
\vskip -1.8  cm
\centerline { \includegraphics[width=.95  \textwidth, angle=0] {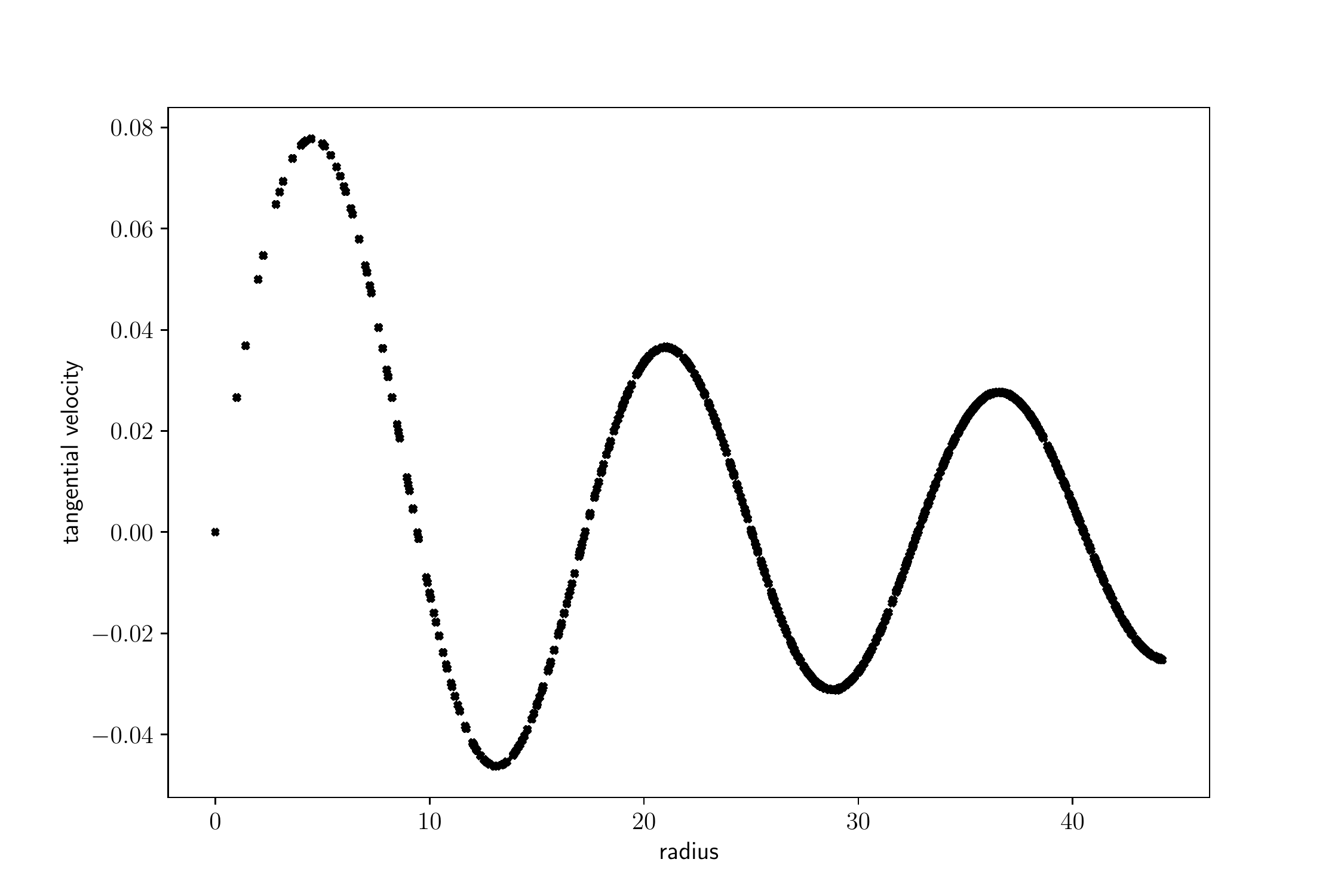}}
\vskip -.5  cm
\caption{Eigenmode   $\Gamma= -1.590  \,$ in unit $ \,  1/\Delta t \,$
for a linear  fluid problem in a disc with homogeneous  anti bounce back boundary condition for a disc of radius
$ \, R = 44.2 $. Tangential velocity profile as a function of the radius. The density and radial velocity are negligible for this mode.   } 
\label{abb-mode-44-densite} \end{figure}

%
\noindent 
In this contribution, 
we have determined the discrete eigenmodes of the lattice Boltzmann
scheme
\moneq \label{lb-modes} 
f_j(x,t+\Delta t) \, \equiv \,  \exp(  \Gamma \, \Delta t) \,\, f_j(x,t)
\monend
for internal or boundary nodes.
For a linear fluid problem, we have determined the eigenmodes (\ref{lb-modes})
with the homogeneous ($ \xi_j \equiv 0 $) anti bounce back scheme  (\ref{fbl02-abb}) at the boundary. 
%
%
We have used the following values of relaxation parameters:
\moneq  \label{s-quartic}  
  s_e = 1.30 \,, \quad    s_x = {{1}\over{{{\sqrt{3}}\over{6}} + {1\over2}}} \,, \quad 
 s_q = {{1}\over{{{\sqrt{3}}\over{3}} + {1\over2}}} \,, \quad  s_d = 1.30  \, . 
 \monend
The first mode is depicted in Figures \ref{abb-mode-34-densite} and \ref{abb-mode-34-vit-radiale}.
An other mode is presented in Figure  \ref{abb-mode-44-densite}.
We have no definitive interpretation of the
extended bounce back and anti  bounce back 
boundary condition (\ref{fbl01-bb}) to  (\ref{fbl02-abb}) in terms of partial
differential equations.
This question is still an open problem at our knowledge.
Nevertheless, some modes are clearly associated to Bessel functions. 
 These modes are invariant by rotation as evident from these figures
 that include data for all points located in the circle.
Recall that this is mainly due to the good choice of the quartic parameters (\ref{s-quartic})
 as presented in \cite{DL11}.

 \bigskip \bigskip   \noindent {\bf \large    6) \quad  
Analysis of anti bounce back for the linear fluid problem }   

\noindent  
The anti  bounce back  boundary condition  
for the fluid problem is very close to the framework presented 
in Section~4  for the scalar case. 
The boundary iteration still follows a scheme of the form ({\ref{cl-generale}). 
The collision  matrix $ \, C \, $ is no longer given by the relation 
 ({\ref{matrice-collision-thermic}), but by    ({\ref{matrice-collision-fluide}). 
The matrices $ \, U \, $ and $ \, T \, $ characterize the locus of the boundary and
the anti bounce back boundary condition. 
They are still given by the relations ({\ref{UU-cl-abb-chaleur}) and 
({\ref{TT-cl-abb-chaleur})  respectively. 
The matrix $ \, K \, $ defined in (\ref{matrice-K}) 
has a new expression:  
\moneq \label{K-fluide-abb}  
 K=\left( \begin{array}{ccccccccc}
\!\!\!\! {{1}\over{6}} (4+\alpha s_e)  & \!\!\!\! \!\!\!\! 0 &  \!\!\!\!  0  &  \!\!\!\!  {{1}\over{6  \lambda^2}} (1-s_e)
&  \!\!\!\! {{1}\over{2 \lambda^2}} (s_x - 1) &  \!\!\!  \!\!\!\! 0 &  \!\!\!\! 0 &  \!\!\!\! 0 &  \!\!\!\! 0 \!\!\!\! \\ 
\!\!\!\! 0 &  \!\!\!\! \!\!\!\! 0  &  \!\!\!\!  0 &  \!\!\!\! 0 &  \!\!\!\! 0 & 
 \!\!\!\! \!\!\!  {{1}\over{\lambda}} (1 - s_x) &   \!\!\!\! 0 &  \!\!\!\! 0 &  \!\!\!\! 0  \!\!\!\! \\ 
\!\!\!\!   {{\lambda}\over{6}} (4+\alpha s_e)  & \!\!\!\!\!\!\!\!  0 &  \!\!\!\!   0  & 
 \!\!\!\!  {{1}\over{6  \lambda}} (1 - s_e)  &  \!\!\!\!  {{1}\over{2 \, \lambda}} (s_x - 1)
&  \!\!\!\! \!\!\! 0 &  \!\!\!\! 0 &   \!\!\!\! 0 &  \!\!\!\! 0  \!\!\!\! \\ 
\!\!\!\!  {{\lambda^2}\over{6}} (4 - 3 \alpha s_e + 2 \beta s_d )
 & \!\!\!\! \!\!\!\!  0 &  \!\!\!\!  0  &   \!\!\!\! {{1}\over{2}} (s_e + 1)
 &  \!\!\!\!    {{1}\over{2}} (1 - s_x) &  \!\!\! \!\!\!\! 0 &    \!\!\!\! 0 &  \!\!\!\! 0 & 
 \!\!\!\! {{1}\over{3 \lambda^2}} (1 - s_d) \!\!\!\! \\ \vspace {-.5cm} \\ 
\!\!\!\!   {{\lambda^2}\over{18}} (-4 + \alpha s_e + 2 \beta s_d )
 &  \!\!\!\! \!\!\!\! 0  &  \!\!\!\! 0 &  \!\!\!\!  {{1}\over{18}} (1 - s_e) &  \!\!\!\!  {{1}\over{2}} (1 + s_x)
 &  \!\!\!\! \!\!\! 0  &  \!\!\!\!  0  &  \!\!\!\! 0 &  \!\!\!\! {{1}\over{9 \lambda^2}}  (1 - s_d) \!\!\!\! \\ 
\!\!\!\! 0 & \!\!\!\!\!\!\!\! 0 &  \!\!\!\! 0 & \!\!\!\!  0 &  \!\!\!\! 0 & \!\!\!  \!\!\!\!1 &
\!\!\!\!  0 &  \!\!\!\!0 &  \!\!\!\!0 \!\!\!\! \\ 
\!\!\!\! 0 & \!\!\!\! \!\!\!\!  \lambda^2 s_q   &  \!\!\!\! 0 &  \!\!\!\! 0 &  \!\!\!\! 0 & 
\!\!\! \!\!\!\! \lambda (1 - s_x) & \!\!\!\!s_q & \!\!\!\!0 & \!\!\!\!0 \!\!\!\! \\ 
\!\!\!\!   {{\lambda^3}\over{3}}  (\alpha s_e + \beta s_d ) & \!\!\!\! \!\!\!\! 0 &  
\!\!\!\! \lambda^2 s_q   &
 \!\!\!\!  {{\lambda}\over{3}} (1-s_e ) &  \!\!\!\!  \!\!\!\!  \!\!\!\! \lambda (1-s_x) &  \!\!\! \!\!\!\! 0 & 
 \!\!\!\! 0 &  \!\!\!\! s_q &  \!\!\!\!  {{1}\over{3 \lambda}} (1-s_d) \!\!\!\! \\  \vspace {-.5cm} \\ 
\!\!\!\!   {{\lambda^4}\over{3}}  (\alpha s_e - 2 \beta s_d ) &  \!\!\!\!\!\!\!\!  0 &  \!\!\!\!  0 &
 \!\!\!\!  {{\lambda^2}\over{3}} (1-s_e) &   \!\!\!\! \lambda^2 (1-s_x) &  \!\!\! \!\!\!\! 0 &  \!\!\!\! 0 & 
 \!\!\!\! 0 &  \!\!\!\! {{1}\over{3}} (1 + 2 s_d) \!\!\!\!
 \end{array} \right) \!  .  \monend

\monitem 
The matrix $ \, K \, $ presented  in (\ref{K-fluide-abb})   is   singular. 
The associated kernel  is  of dimension 2.  
It is generated by the following two vectors $ \, \kappa_x \, $ and 
$ \, \kappa_y\, $ in the space of moments:
\moneq \label{ker-K-fluide-abb}  
\kappa_x =  \big( 0 ,\, 1 ,\, 0 ,\, 0 ,\, 0 ,\, 0 ,\, -\lambda^2 ,\, 0 ,\, 0 \big)^{\rm  t}  \,, \,\, 
\kappa_y =  \big( 0 ,\, 0 ,\, 1 ,\, 0 ,\, 0 ,\, 0 ,\, 0 ,\, -\lambda^2 ,\, 0  \big)^{\rm t}  \, . 
\monend 
This means first that when solving a generic linear system of the type 
\moneq \label{systeme-modele}  
K \, m = g \equiv
 \big( g_\rho ,\, g_{jx} ,\, g_{jy} ,\, g_\varepsilon ,\, g_{\varphi x} ,\, g_{\varphi y} ,\, g_{q x}  ,\, g_{q x} ,\, g_h \big)^{\rm  t}  \,, 
\monend
two compatibility  relations  have to be satisfied by    
the right hand side $ \, g $: 
\moneq \label{compatibilite-fluide-abb}  
\lambda \, g_\rho - g_{jy} = 0  \,, \,\, 
\lambda \, g_{jx} + (s_x-1) \, g_{\varphi y} = 0  \, .  
\monend
Secondly, it is always possible to add to any solution of the model system    
(\ref{systeme-modele}) any combination of the type $ \, j_x \, \kappa_x +  j_y \, \kappa_y $. 
In other terms, two components  $ \, j_x \, $ and $ \, j_y \, $ of momentum   remain  undetermined
by the anti bounce back boundary condition.  They have to be evaluated in practice by
the interaction with the  other vertices through the numerical scheme.

\monitem  
Following the procedure presented  previously  for the thermal case, we can try to solve 
the boundary condition (\ref{abb-systeme-lineaire}) at various orders. 
At order zero,  no relevant information is given by the  compatibility  conditions
 (\ref{compatibilite-fluide-abb}). Then the solution $ \, m_0 \, $
at order zero depends on two parameters identified as the momentum
$ \, (j_x ,\, j_y ) \, $ in the first cell. With a given pressure $ \, p_0 \, $ 
or a  given density   $ \, \rho_0 \equiv {{p_0}\over{c_0^2}} \, $   on the boundary, 
we have 
\moneq \label{m0-fluide-abb}  
m_0=\left(\rho_0, \, j_x, \,j_y, \,\alpha  \, \lambda^2\, \rho_0  ,\, 0,\, 0,\, 
- \lambda^2 \, j_x ,\, - \lambda^2 \, j_y , \,
 \beta \, \,  \lambda^4 \rho_0  \right )^{\rm t}  \, . 
\monend
At  order one, the first compatibility relation  is a linear combination 
 of the first order equivalent partial differential  equations. 
The other compatibility relation gives a 
  non trivial differential condition  on the boundary : 
\moneq \label{compatibilite-fluide-abb-final}
\partial_x J_y + \partial_y J_x  = 0 \, .  
\monend
%
The anti bounce back induces a hidden
additional boundary condition. 
Observe that the differential condition (\ref{compatibilite-fluide-abb-final})  is not satisfied for a Poiseuille flow.
Nevertheless, if this relation 
is satisfied, it is possible to expand the resulting density in the first cell:
\moneq \label{dentite-ordre-un-fluide-abb}   
\rho \Big( x_1, {{\Delta x}\over{2}} \Big)  = \rho_0 \big( x_1 \big) +   {{1}\over{2}} \, \Delta x \, \, \partial_y \rho(x_1) + \Delta t \,  \big(  
a_{jx} \,\, \partial_x J_x + a_{jy} \,\, \partial_y J_y \big) + {\rm O}(\Delta x^2)  \, . 
\monend
with
\moneqstar \left \{ \begin{array} {l} \displaystyle 
a_{jx}  = {{1}\over{ 8 \, (s_e \, s_x + s_e \, s_d + s_x \, s_d) + 2 \, (s_e \, s_x + 2 \, s_e \, s_d) \, \alpha + 2 \, (  s_e \, s_d  - s_x \, s_d) \, \beta }} \, 
\Big(  6 \, s_e \, s_x + 5 \, s_e \, s_d \\ \displaystyle \qquad  \qquad 
+ \, 7 \, s_x \, s_d - 6 \, s_e + 2 \, s_x - 8 \, s_d + 2 \, (s_e \, s_x + 2 \, s_e \, s_d - s_x - 2 \, s_d) \, \alpha
\\ \displaystyle \qquad  \qquad  + 2 \, (s_e \, s_d  - s_x \, s_d  - s_e + s_x) \, \beta \Big)
\\ \displaystyle   \vspace{-4 mm} \\ \displaystyle   
a_{jy}  = {{1}\over{ 8 \, (s_e \, s_x + s_e \, s_d + s_x \, s_d) + 2 \, (s_e \, s_x + 2 \, s_e \, s_d) \, \alpha + 2 \, (  s_e \, s_d  - s_x \, s_d) \, \beta }} \, 
\Big( 2 \, s_e \, s_x + 5 \, s_e \, s_d\\ \displaystyle \qquad  \qquad 
- \, s_x \, s_d + 2 \, s_e + 2 \, s_x + 8 \, s_d  + 2 \, (s_e \, s_x + 2 \, s_e \, s_d - s_x - 2 \, s_d) \, \alpha 
\\ \displaystyle \qquad  \qquad  + 2 \, (s_e \, s_d  - s_x \, s_d  - s_e + s_x) \, \beta \Big) \, . 
\end{array} \right . \monendstar 
This result indicates that the anti bounce back boundary condition gives not
satisfactory results for the field near the boundary. 
 The extra term $ \, a_{jx} \,\, \partial_x J_x + a_{jy} \,\, \partial_y J_y \, $
in the relation (\ref{dentite-ordre-un-fluide-abb}) is a discrepancy if we compare this relation with a
common Taylor expansion of the density in the normal direction as found in the scalar case
   (\ref{abb-thermique-ordre-un}).
The relations  (\ref{compatibilite-fluide-abb-final}) and 
(\ref{dentite-ordre-un-fluide-abb}) put  in evidence
quantitatively various  defects 
of the anti bounce back boundary condition as proposed at the relations 
(\ref{abb-chaleur}). 

\newpage 
\monitem  Numerical experiments for a linear  Poiseuille flow.

\noindent
We have considered a two-dimensional vertical channel  with wall boundaries at the left and at the right.
At the bottom of the channel a given $ \, +P_0 \, $ pressure is imposed through anti bounce back 
and  $ \, -P_0 \, $ is imposed at the top (see Figure \ref{test-poiseuille}). 


\begin{figure}    [H]  \centering 
\vskip -1.8  cm
\centerline { \includegraphics[width=.45  \textwidth, angle=0] {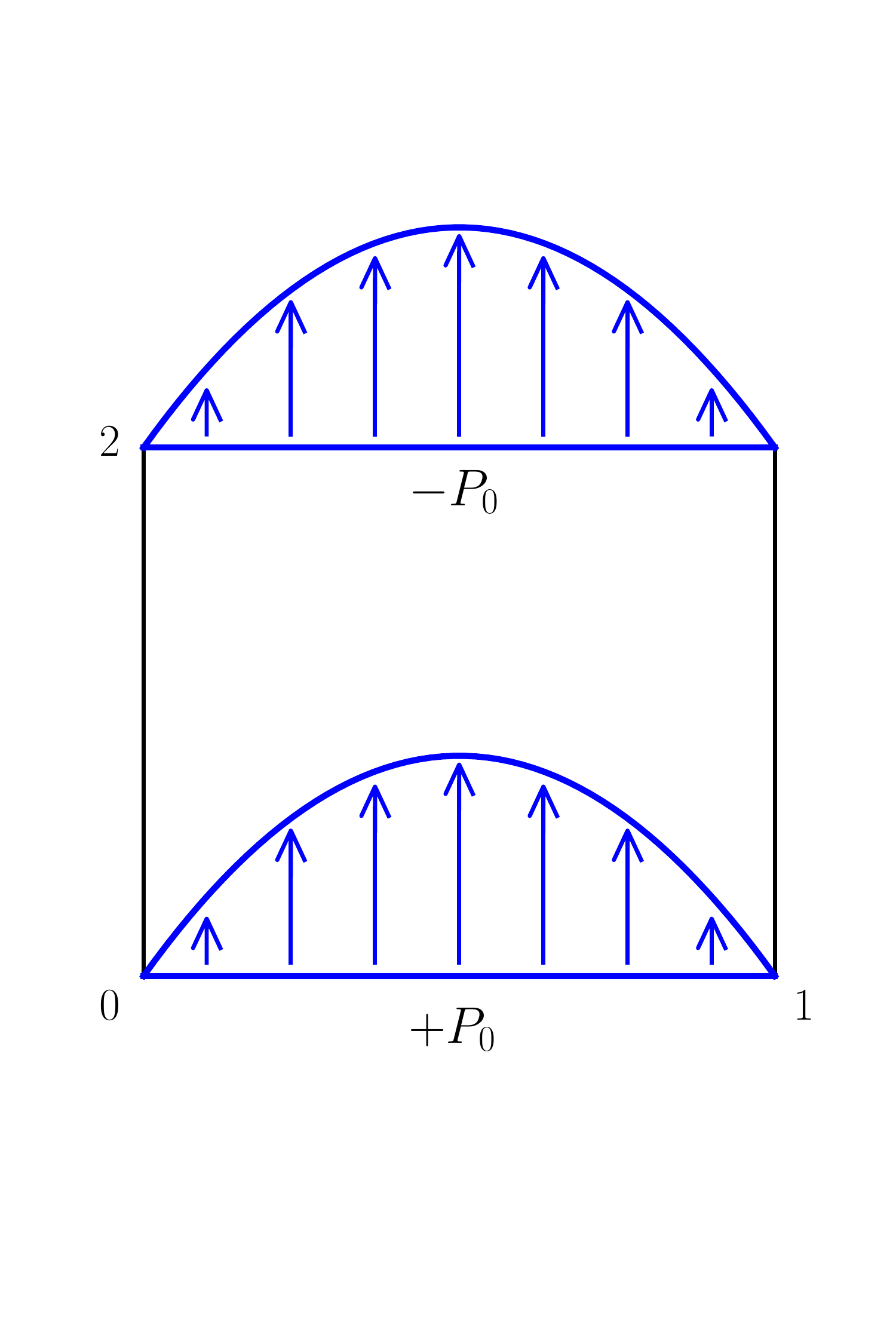}}
\vskip -1.5 cm
\caption{Linear vertical Poiseuille flow. The solid boundaries are at the left and at the right.
  A given pressure of $ \, P_0 = 0.001 \, $ and $ \, - P_0 \, $ are given at the bottom and the top of the channel. 
  A parabolic profile is generated by the pressure gradient. } 
\label{test-poiseuille} \end{figure}

\noindent
We have emphasized in Figures \ref{abb-poiseuille-XY}, \ref{abb-poiseuille-pression}, \ref{abb-poiseuille-vitesse-x} and \ref{abb-poiseuille-vitesse-y}
the results for  three meshes with $ 20 \times 40 $, $ 40 \times 80 \, $ and $ \, 80 \times 160 \, $ cells.
In Figure \ref{abb-poiseuille-XY}, the field is represented in the first two layers near the inflow boundary.
Due to the classical Taylor expansion \cite{LL00}
\moneqstar
\varphi_y = - {{\lambda \, \Delta x}\over{3 \, s_x}} \,  \big( \partial_x J_y + \partial_y J_x \big) + {\rm O}(\Delta x^2) \,,
\monendstar
this field is an excellent indicator of the treatment of the hidden boundary condition (\ref{compatibilite-fluide-abb-final}).
The Figure \ref{abb-poiseuille-XY} confirms that for the first cell, the hidden condition is effectively taken into
account, except may be  in the corners. 

\noindent
In  consequence, the Poiseuille flow cannot be correctly satisfied near the fluid boundary.  
An important error  of  18.2\%  for the pressure in the first cell is observed 
for the three meshes (see Fig. \ref{abb-poiseuille-pression}). Nevertheless, if we fit these pressure results 
by linear curves taking only half the number of points in the middle of the channel, this error is reduced to  6.5\% 
(see Fig.~\ref{abb-poiseuille-pression}).
The tangential velocity is supposed to be identically null in all the channel. 
In the middle, the maximal  tangential velocity is always less than 0.01\% of the maximum normal velocity.
But in the first cells (Fig. \ref{abb-poiseuille-vitesse-x}), the discrepancies are important
and can reach 17\% of the maximal velocity. 
%
%

\begin{figure}    [H]  \centering 
\vskip -1.8  cm
\centerline { \includegraphics[width=.95  \textwidth, angle=0] {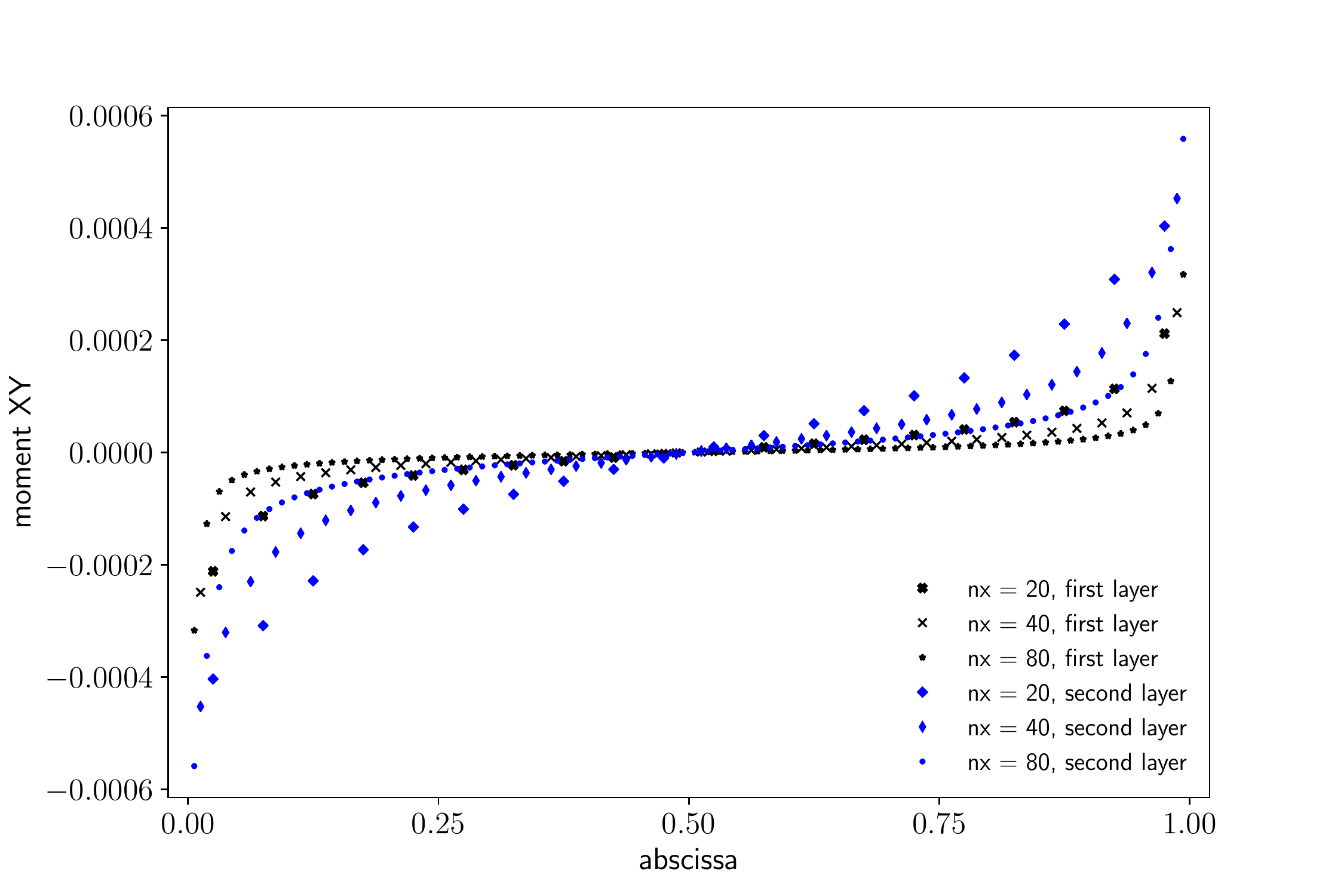}}
\vskip -.5 cm
\caption{Anti bounce back boundary condition for linear Poiseuille flow. Moment  
field $\, \varphi_y $
in the first layers of the boundary. The hidden boundary condition $ \, \partial_x J_y + \partial_y J_x = 0 \, $
enforces the constraint $ \,  XY \equiv \varphi_y = 0 $. } 
\label{abb-poiseuille-XY} \end{figure}

\begin{figure}    [H]  \centering 
\vskip -1.8  cm
\centerline { \includegraphics[width=.95  \textwidth, angle=0] {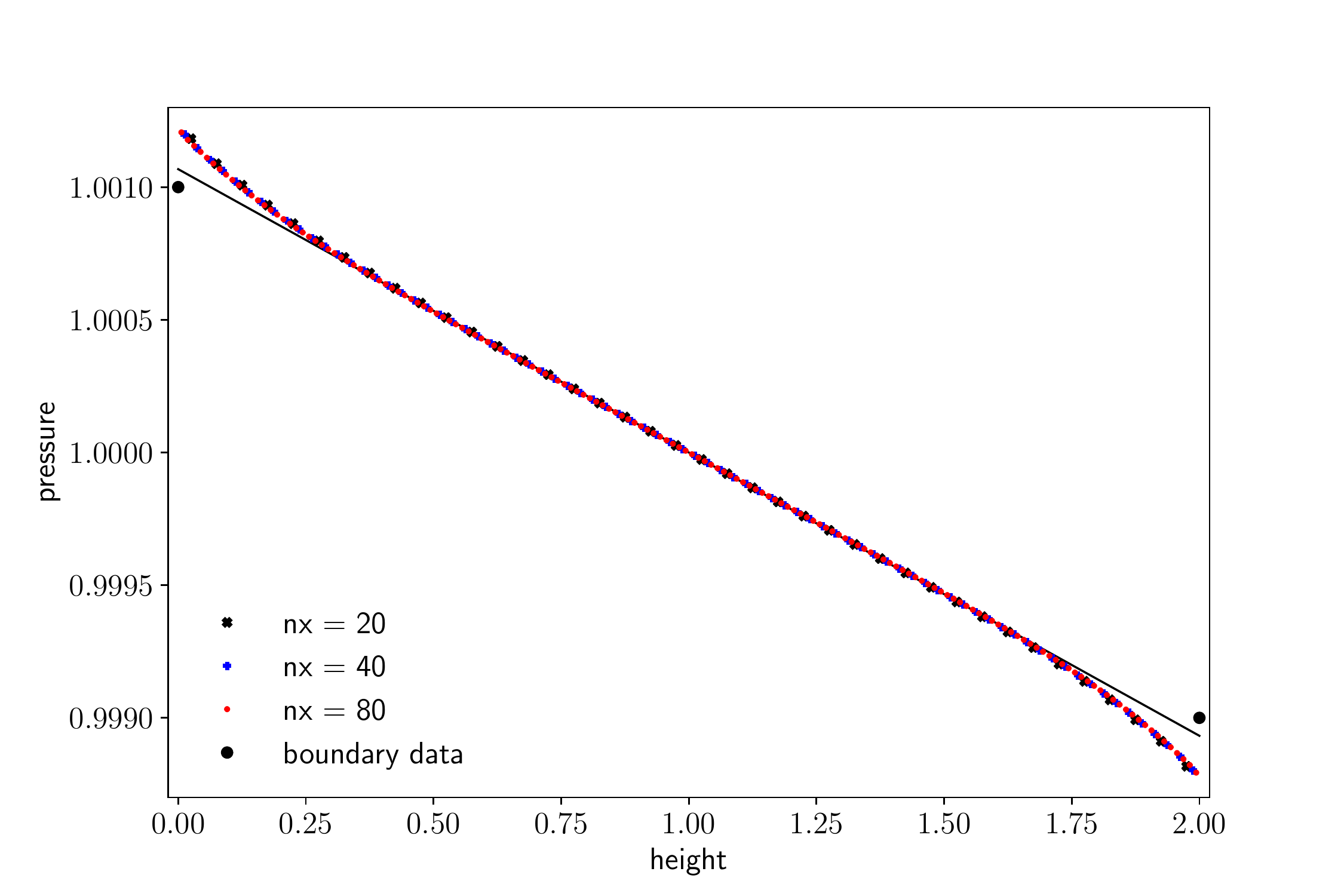}}
\vskip -.5 cm
\caption{Anti bounce back boundary condition for linear Poiseuille flow.
Pressure field in the middle of the channel (after a simple interpolation due to an even number of meshes) 
for three meshes with $ 20 \times 40 $, $ 40 \times 80 \, $ and $ \, 80 \times 160 \, $ cells.  
 A  substantial error does not vanish  as the mesh size tends to zero. } 
\label{abb-poiseuille-pression} \end{figure}

\begin{figure}    [H]  \centering 
\vskip -1.8  cm
\centerline { \includegraphics[width=.95  \textwidth, angle=0] {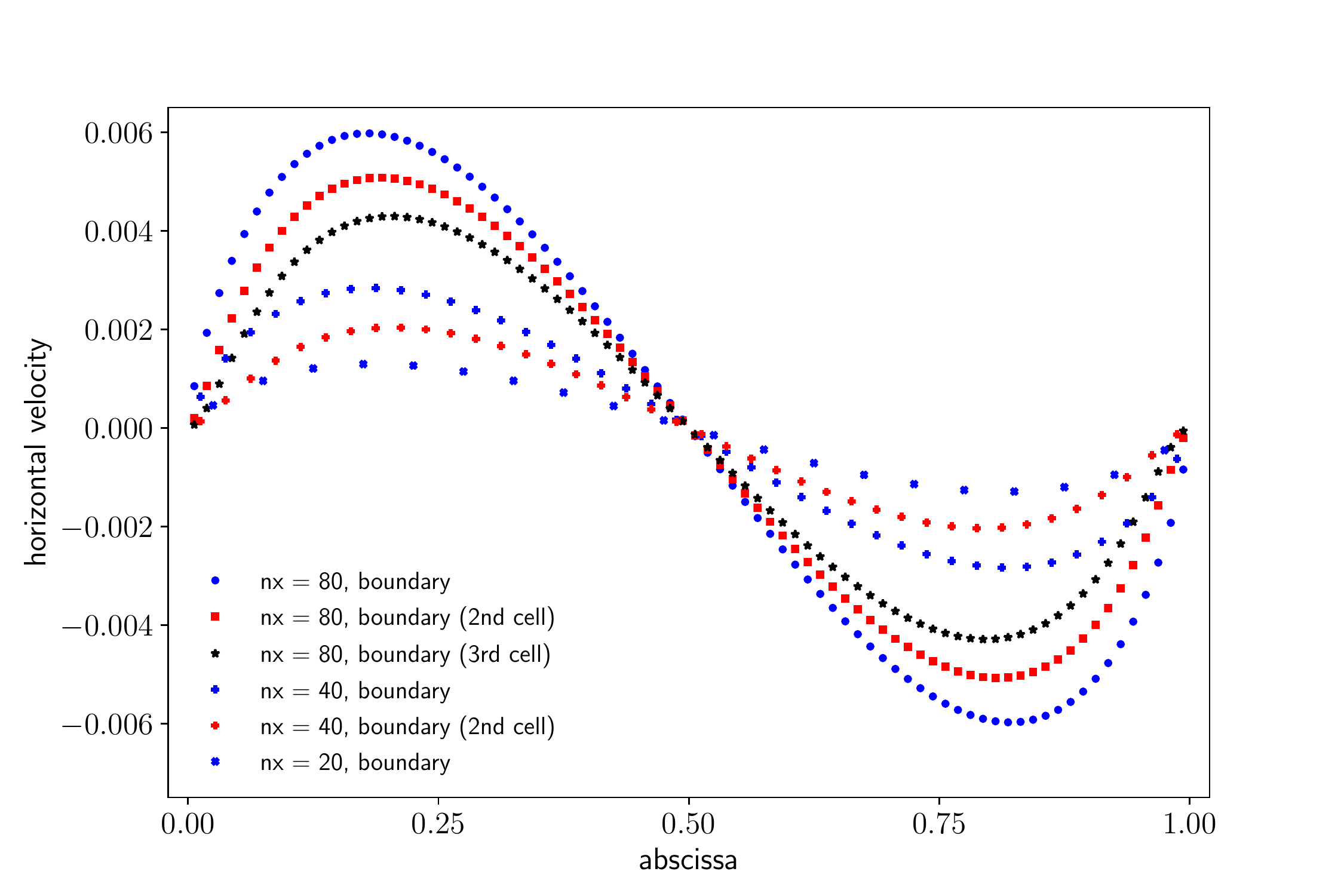}}
\vskip -.5 cm
\caption{Anti bounce back boundary condition for linear Poiseuille flow.
Horizontal velocity field for three meshes. The maximum value is around 17\% of the vertical velocity.   }  
\label{abb-poiseuille-vitesse-x} \end{figure}

\begin{figure}    [H]  \centering 
\centerline { \includegraphics[width=.95  \textwidth, angle=0] {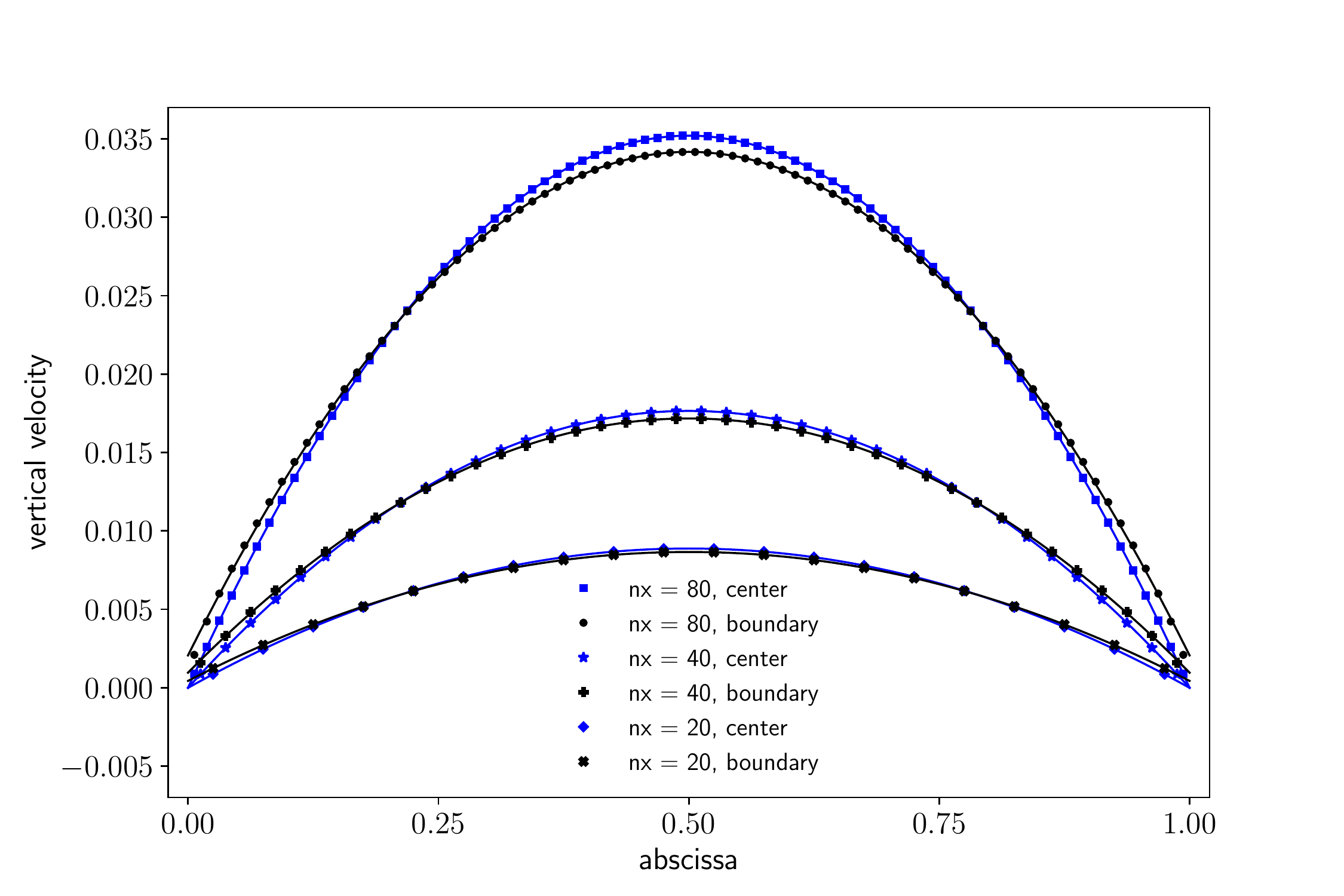}}
\vskip -.5 cm
\caption{Anti bounce back boundary condition for linear Poiseuille flow.
Vertical velocity field for three meshes. The parabolic profile is not  completely recovered, in particular
in the cells near the fluid boundary. } 
\label{abb-poiseuille-vitesse-y} \end{figure}
%

%
\noindent
In Figure \ref{abb-poiseuille-vitesse-y}, the axial velocity is presented for the three meshes.
We keep  the same value of $ \, s_x \, $ and in consequence the kinematic viscosity of the fluid \cite{Du08, LL00} 
\moneqstar
\nu = {{\lambda}\over{3}} \, \Delta x  \, \Big( {{1}\over{s_x}} - {1\over2} \Big) \, 
\monendstar
is reduced at each mesh refinement by a factor of 2
and the normal velocity is multiplied in the same proportions.
We have fitted with least squares the vertical velocity by a parabola in the middle of the
channel and in the first layer (see Figure \ref{abb-poiseuille-vitesse-y}). 
The null velocity along this fitted parabola defines a numerical boundary that
is compared to its true geometric location. The results are presented in  Table~\ref{table-abb}. 
In the center of the channel,
we obtain a typical  error of $\, 10^{-3} \, $ between the numerical and geometrical  boundaries,
measured in cell units: 
the axial velocity is of very good quality in the center of the channel.
In the first cell, 
the gap between  numerical and physical boundaries represents
a notable fraction of one unit cell (see Table~\ref{table-abb}). 

%
\begin{table}  [H]     \centering
 \centerline { \begin{tabular}{|c|c|c|c|}    \hline 
 mesh &  $ 20 \times 40 $ & $ 40 \times 80  $ & $   80 \times 160  $ \\   \hline 
 center & $ 5.29  \, 10^{-3}   $ & $  3.98 \, 10^{-3}  $ & $  4.61 \, 10^{-3}  $  \\   \hline 
 bottom   & $  0.261  $  & $  0.582 $ &  $ 1.26 $  \\   \hline  
\end{tabular} }  
 \caption{Anti bounce back boundary condition for linear Poiseuille flow.
 Distance between the numerical position of the boundary from its theoretical location. } 
\label{table-abb} \end{table}
%

 \bigskip \bigskip   \noindent {\bf \large    7) \quad  
Mixing bounce back and anti bounce back boundary conditions}   

\noindent  
In order to overcome the difficulties with the  anti bounce back boundary condition 
for the fluid problem, we have adapted a mixing of  bounce back and anti bounce back  
first proposed in \cite{DL15} for the thermal Navier-Stokes equations. 
We suppose that all the fluid information, {\it id est} density and momentum, 
is given on the boundary. We force this relation by taking bounce back boundary condition 
for the particles  coming from the left  and from the right:
 $ \, f_5^{\rm eq} -  f_7^{\rm eq} = {{1}\over{6 \, \lambda}} \, ( J_x + J_y ) \, $ and 
$ \,   f_6^{\rm eq} -  f_8^{\rm eq} = {{1}\over{6 \, \lambda}} \, ( -J_x + J_y ) $.  
We keep the  anti bounce back 
$ \,   f_2^{\rm eq} +  f_4^{\rm eq} = {{4 -\alpha-2 \, \beta}\over{18}} \, \rho \, $ 
for the particles coming from the bottom, as illustrated in Figure~\ref{mixing-bb-abb}. 
%
\begin{figure}    [H]  \centering 
  \centerline { \includegraphics[width=.35 \textwidth, angle=0] {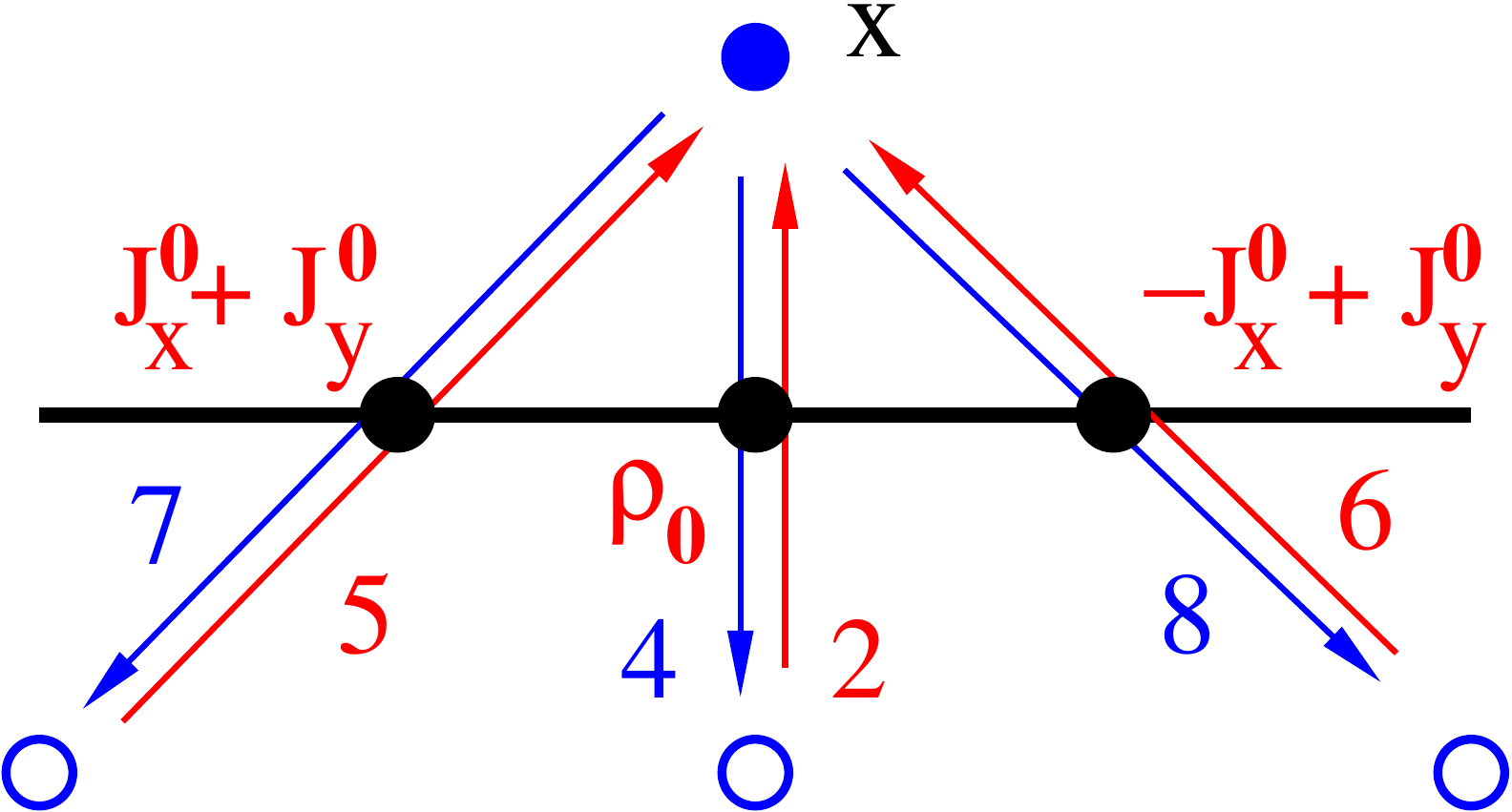}} 
\caption{Mixed bounce back and anti bounce back boundary condition to enforce 
density and momentum values on the boundary.    } 
\label{mixing-bb-abb} \end{figure}
%
Finally, if density $\, \rho_0 \, $ and the two components $ \, J_x^0  \, $ 
and $ \, J_y^0  \, $  of the velocity  are supposed  given on the boundary, 
the mixed bounce back and anti bounce back boundary condition can be written as: 
\moneq \label{mixed-bb-abb-fluide}  \left \{ \begin{array} {l} \displaystyle 
f_{5}(x, \, t+ \Delta t) =   f^*_{7}(x,\, t) \,+\,   {{1}\over{6 \, \lambda}} \, 
\big( J_x^0 + J_y^0 \big) \Big( x_1 - {{\Delta x}\over{2}} , \, t  \Big)    \,, 
\\  \vspace{-4 mm} \\ \displaystyle
f_{2}(x,  \, t+ \Delta t) =    -    f^*_{4}(x,\, t) \,+\, {{4 -\alpha-2 \, \beta}\over{18}} \, 
\rho_0 \big(x_1 , \, t \big)  \,,  
\\  \vspace{-4 mm} \\ \displaystyle 
\displaystyle f_{6}(x,  \, t+ \Delta t)  =  f^*_{8}(x ,\, t) \,+\,  {{1}\over{6 \, \lambda}} \, 
\big( -J_x^0 + J_y^0 \big) \Big( x_1 + {{\Delta x}\over{2}} , \, t  \Big)   \, . 
\end{array} \right.  \monend

\monitem 
The analysis of the discrete mixed condition (\ref{mixed-bb-abb-fluide})
follows the methodology presented above. 
The boundary iteration still follows the general framework ({\ref{cl-generale}). 
The matrices $ \, U \, $ and $ \, T \, $ characterize the locus of the boundary and
the mixed bounce back and anti bounce back boundary condition. 
The matrix $ \, U \, $ is still given by  ({\ref{UU-cl-abb-chaleur}).
On the contrary, the interaction matrix $ \, T \, $ is no longer given by the expression  
({\ref{TT-cl-abb-chaleur}), but by the relation  
\moneq \label{TT-cl-mixed-bbabb-fluide} 
 T = \left( 
\begin{array}{ccccccccc}
0 & 0 & 0 & 0 & 0 & 0 & 0 & 0 & 0 \\ 
0 & 0 & 0 & 0 & 0 & 0 & 0 & 0 & 0 \\ 
0 & 0 & 0 & 0 & \! -  1 & 0 & 0 & 0 & 0 \\ 
0 & 0 & 0 & 0 & 0 & 0 & 0 & 0 & 0 \\ 
0 & 0 & 0 & 0 & 0 & 0 & 0 & 0 & 0 \\ 
0 & 0 & 0 & 0 & 0 & 0 & 0 & \!  + 1 & 0 \\ 
0 & 0 & 0 & 0 & 0 & 0 & 0 & 0 & \! +  1 \\ 
0 & 0 & 0 & 0 & 0 & 0 & 0 & 0 & 0  \\ 
0 & 0 & 0 & 0 & 0 & 0 & 0 & 0 & 0 
\end{array} \right) \, . \monend
The matrix $ \, K \, $ introduced  in (\ref{matrice-K}) 
takes the form:  
\moneq \label{K-fluide-mixed-bbabb}
K=\left( \begin{array}{ccccccccc}
\!\!\!\! {{4-\alpha s_e-2 \beta s_d}\over{18}}  & \!\!\! \!\!\!\! 0 &  \!\!\!\! \!\!\!\!   {{2 - s_q }\over{3 \lambda}}
  &  \!\!\!   \!\!  {{s_e-1}\over{18  \lambda^2}} &  \!\!\!\!  {{s_x-1}\over{2  \lambda^2}} 
&  \!\!\!\! \!\!\! 0 &  \!\!\!  \!\!\!\! 0 &  \!\!\!\! {{1 - s_q}\over{3 \, \lambda^3}} 
 &  \!\!\!\! {{s_d -1 }\over{9 \lambda^4}}   \\ 
\!\!\!\! 0 & \!\!\!  \!\!\!\!   {{2 - s_q }\over{3}}   & \!\!\!\!   \!\!\!\!  0 & \!\!\!   \!\!  0 &  \!\!\!\! 0 & 
 \!\!\!\! \!\!\!  0 &   \!\!\!\!  {{1 - s_q }\over{3 \lambda^2}}  &  \!\!\!\! 0 &  \!\!\!\! 0  \\ 
\!\!\!\!    {{\lambda  (4 - \alpha s_e - 2 \beta s_d )}\over{18}} & \!\!\! \!\!\!\!  0 &  
 \!\!\!\! \!\!\!\!    {{2 - s_q }\over{3}}   & 
 \!\!\!\!\!  {{s_e -1  }\over{18  \lambda}}  &  \!\!\!\!  {{s_x - 1}\over{2  \lambda}}  
&  \!\!\!\! \!\!\! 0 &  \!\!\!\! 0 &   \!\!\!\!  {{1 - s_q }\over{3 \lambda^2}}  &  
\!\!\!\!   {{s_d -1 }\over{9 \lambda^3}}   \\ 
\!\!\!\!  {{\lambda^2  (-4 - 17 \alpha s_e + 2 \beta s_d )}\over{18}} 
 & \!\!\! \!\!\!\!  0 &  \!\!\!\!  \!\!\!\!   {{\lambda (4 - 2 s_q) }\over{3}} & \!\!\!  
 \!\! {{1 + 17 s_e }\over{18}} 
 &  \!\!\!\!    {{1 - s_x}\over{2}} &  \!\!\! \!\!\!\! 0 &    \!\!\!\! 0 &  
\!\!\!\!  {{2 (1 - s_q) }\over{3 \lambda}}  &  \!\!\!\! {{1 - s_d }\over{3 \lambda^2}} \\ 
\!\!\!\!   {{\lambda^2  (-4 + \alpha s_e + 2 \beta s_d )}\over{18}} 
 & \!\!\! \!\!\!\! 0  &  \!\!\!\! \!\!\!\! 0  &  \!\!\!\!\!   {{1 - s_e }\over{18}} &  \!\!\!\!  {{1 + s_x}\over{2}} 
 &  \!\!\!\! \!\!\! 0  &  \!\!\!\! 0  &  \!\!\!\! 0 &  \!\!\!\! {{1 - s_d}\over{9 \lambda^2}}  \\ 
\!\!\!\! 0 & \!\!\!\!  \!\!\!\!  {{\lambda (2-s_q )}\over{3}} 
 &   \!\!\!\! \!\!\!\! 0 & \!\!\! \!\!  0 &  \!\!\!\! 0 & \!\!\!  \!\!\!\! s_x &
  \!\!\!\!   {{1-s_q }\over{3 \lambda}}   &  \!\!\!\!0 &  \!\!\!\!0 \\ 
\!\!\!\! 0 & \!\!\! \!\!\!\!  {{2 \lambda^2 (1+s_q) }\over3}   &  \!\!\!\!  \!\!\!\! 0 & \!\!\!  \!\!  0 &  
\!\!\!\! 0 & 
\!\!\! \!\!\!\! 0  & \!\!\!\!   {{1 + 2 s_q  }\over3}  & \!\!\!\!0 & \!\!\!\!0 \\ 
\!\!\!\!   {{\lambda^3  (-4 + \alpha s_e + 2 \beta s_d )}\over{9}}  & \!\!\! \!\!\!\!  0 &  
 \!\!\!\! \!\!\!\!     {{2 \lambda^2 (1+s_q) }\over3} &
 \!\!\! \!\!  {{\lambda (1-s_e )}\over{9}} &  \!\!\!\! \lambda (1 \!\! - \!\! s_x\!) &  \!\!\! \!\!\!\! 0 & 
 \!\!\!\! 0 &  \!\!\!\!  {{1 + 2 s_q  }\over3} &  \!\!\!\!  {{2(1-s_d) }\over{9 \lambda}} \\ 
\!\!\!\!   {{\lambda^4  (-4 + \alpha s_e - 7 \beta s_d )}\over{9}}  &  \!\!\! \!\!\!\!  0 &  
 \!\!\!\! \!\!\!\!   {{\lambda^3 (2-s_q )}\over{3}}  &
\!\!\!  \!\!   {{\lambda^2 (1-s_e)}\over{3}} &   \!\!\!\! \lambda^2 (1 \!\! - \!\! s_x \!) &  \!\!\!\! 0 &
 \!\!\!\! 0  & \!\!\! \!\!\!\!  {{\lambda (1-s_q)}\over{3}}  &  \!\!\!\! {{2 + 7 s_d}\over{9}} 
\end{array} \right) \, .  \monend
The matrix 
defined in (\ref{K-fluide-mixed-bbabb}) 
is   singular. Its   kernel  is of dimension 1. It is  generated by the following vector $\, \kappa $:
\moneqstar  
\kappa = \begin{pmatrix}  \displaystyle 
{{(3 \, s_e \, s_d+2 \, s_e \, s_x+s_d \, s_x) \, s_q }
\over{\lambda \,  s_d\, s_e \, s_x \, (\alpha+2 \, \beta-4) }} \\  \vspace{-4 mm} \\ \displaystyle
0 \\ 1 \\ \displaystyle {{ \lambda \, (3 \, \alpha \, s_e \, s_d+2 \, \alpha \, s_e \, s_x-2 \, \beta \, s_d \, s_x
+4 \, s_d \, s_x) \, s_q }\over{ s_d \, s_e \, s_x \, (\alpha+2 \, \beta-4)}}  \\   \vspace{-4 mm} \\ \displaystyle
{{\lambda \, s_q  }\over{3 \, s_x}} \\ 
0 \\ 0  \\ \displaystyle 
-2 \, \lambda^2 \\ \vspace{-4 mm} \\ \displaystyle   -\lambda^3 \, {{ \alpha \, s_e \, s_x-3 \, 
\beta \, s_e \, s_d-\beta \, s_d \, s_x-4 \, s_e \, s_x) \, s_q } 
\over{ s_d \, s_e \, s_x \, (\alpha+2 \, \beta-4)}} \\
\end{pmatrix} \, .  \monendstar
We did not find any  simple physical interpretation of the kernel $ \, \RR \kappa \, $ of the matrix $ \, K \, $ 
defined by  (\ref{K-fluide-mixed-bbabb}) in this case. 
When solving the generic linear system (\ref{systeme-modele}), 
the right hand side $ \, g \, $ must satisfy 
\moneqstar  
\lambda \, g_\rho - g_{jy} = 0  \, .  
\monendstar
Moreover, the solution is defined up to a multiple of the vector  $ \, \kappa \, $ 
presented above.  

\monitem 
At order zero, the compatibility condition does not give any relevant information.
If we suppose that the component of $\, m_0 \, $ along the eigenvector $ \, \kappa \, $ 
is reduced to zero, we have   for a given density  $ \, \rho_0 \, $  
and a given momentum $ \, ( J_x^0  , \, J_y^0 ) \, $ 
on the boundary: 
\moneqstar  
 m_0=\left(\rho_0, \, J_x^0 , \,J_y^0 , \,\alpha  \, \lambda^2\, \rho_0  ,\, 0,\, 0,\, 
- \lambda^2 \, J_x^0   ,\, - \lambda^2 \, J_y^0 , \,
\beta \, \,  \lambda^4 \rho_0  \right )^{\rm t}  \, . 
\monendstar
At   order one  the compatibility relations is a  linear combination 
 of the first order equivalent partial differential  equations  
and no constraint is added by this relation. 
The resulting moments in the first cell $ \, x \equiv  \big(  x_1,\,  {{\Delta x}\over{2}} \big) \, $ 
at first order can be expressed  after some formal calculus:

\moneq   \label{bbabb-champs-premiere-maille} 
\left\{ \begin{array}{l} \displaystyle 
  \rho \Big( x_1,\,  {{\Delta x}\over{2}} \Big)  = \rho_0 \big( x_1 \big) + {{1}\over{2}} \, \Delta x \, \partial_y \rho \big( x_1 ,\, 0 \big)
  +  {{ \Delta x}\over{s_x \, s_e \, s_d  \, (\alpha + 2\,\beta-4)}} \, \Big( a_\rho  \, \, \partial_y \rho \big( x_1 ,\, 0 \big)
\\ \vspace{-4 mm} \\ \displaystyle \qquad \qquad \qquad \qquad \qquad 
  + \, b_\rho \,\,   \partial_x J_x  \big( x_1 ,\, 0 \big) + c_\rho \,\,   \partial_y J_y  \big( x_1 ,\, 0 \big)  \Big)  +  {\rm O} (\Delta x^2) \,, 
\\ \vspace{-4 mm} \\ \displaystyle 
j_x  \Big( x_1, \, {{\Delta x}\over{2}} \Big) 
= J_x^0 \big( x_1 \big)   + {1\over2} \, \Delta x \, \, \partial_y J_x \big( x_1 ,\, 0 \big) +  {\rm O} (\Delta x^2) \,, 
\end{array} \right.  \monend
with 
\moneq   \label{bbabb-coefficients-champs} 
\left\{ \begin{array}{l} \displaystyle  
a_\rho = {1\over6} \, \Big(- 4 \, (2 \, s_e \, s_x + 3 \, s_e \, s_d + s_x \, s_d) 
\\ \vspace{-4 mm} \\ \displaystyle \qquad \qquad \qquad 
+ \,(2 \, s_e \, s_x \, s_q + 3 \, s_e \, s_q \, s_d + s_x \, s_q \, s_d  - 6 \, s_e \, s_x - 9 \, s_e \, s_d - 3 \, s_x \, s_d) \, \alpha
\\ \vspace{-4 mm} \\ \displaystyle \qquad \qquad \qquad 
+ \,(2 \, s_e \, s_x \, s_q +  3 \, s_e \, s_q \, s_d + s_x \, s_q \, s_d  - 4 \, s_e \, s_x   - 6 \, s_e \, s_d - 2 \, s_x \, s_d) \, \beta  \Big) \, , 
\\ \vspace{-4 mm} \\ \displaystyle 
b_\rho =  2 \, s_e \, (3 \, s_d - 2 \, s_x \, s_d - s_x) + s_x \, s_d \, (s_e-1) \,  \alpha + 2 \, s_e \, s_x \, (s_d-1) \, \beta \, , 
\\ \vspace{-4 mm} \\ \displaystyle 
c_\rho = {1\over2} \, \Big( 4 \, s_e \, s_x \, s_d  - 2 \, s_e \, s_x \, s_q - 3 \, s_e \, s_q \, s_d - s_x \, s_q \, s_d - 4 \, s_e \, s_x - 12 \, s_e \, s_d
\\ \vspace{-4 mm} \\ \displaystyle \qquad \qquad \qquad 
+ 2 \, s_x \, s_d \, (s_e-1) \,  \alpha + 4 \, (s_d-1) \, s_e \, s_x \, \beta \Big) \,  
\end{array} \right.  \monend
and
\moneq   \label{bbabb-jy-premiere-maille} 
j_y  \Big( x_1, \, {{\Delta x}\over{2}} \Big) = J_y^0 \big( x_1 \big) +  {\rm O} (\Delta x^2) \,.
\monend
We are puzzled by this result and more work is needed.  
There is a true discrepancy for the density with the three terms parametrized  by the 
coefficients $ \, a_\rho $,  $ \, b_\rho \, $ and  $ \, c_\rho $.
Moreover, no first order term $ \, {1\over2} \, \Delta x \, \partial_y J_y \, $ 
of a simple  Taylor formula is present for the expansion (\ref{bbabb-jy-premiere-maille}) of the normal momentum.

 \bigskip \bigskip   \noindent {\bf \large    8) \quad  
Giving pressure and tangential velocity on the boundary }   

\noindent  
A natural mathematical question is to know whether a given set of partial differential equations
associated with a given set of boundary conditions conducts to a well posed problem   
in the sense of Hadamard.
For example, the Poisson  equation with Dirichlet or Neumann boundary conditions 
conducts to a well posed problem (see {\it e.g.} \cite{Ev98}).
On the contrary for the same Laplace equations, the Cauchy problem is defined by the fact to impose
the value of the unknown field and the value of the normal derivative on some part of the boundary.
As well known, this Cauchy  problem for the Laplace equation is not correctly posed \cite {Ca58, Do53}. 
For the stationary Stokes problem  
\moneqstar  
 {\rm div } u = 0 \,, \quad -\nu \, \Delta u + \nabla  p = 0  \,, 
\monendstar 
a set of correctly posed boundary conditions is based on the 
velocity vorticity pressure formulation \cite{Du92}:
\moneq \label{stokes-omega-u-p}  
 {\rm div } u = 0 \,,  \quad  
\omega - {\rm curl } \, u = 0 \,,  \quad   
 \nu \, {\rm curl } \, \omega + \nabla p  = 0  \, . 
\monend 
Considering a variational formulation of the Stokes system  (\ref{stokes-omega-u-p}), 
natural boundary conditions 
can be derived and we obtain the following procedure.   
Consider the two following  decompositions $ \, \Gamma_m \cup \Gamma_p \, $
and $ \, \Gamma_t \cup \Gamma_\theta \, $   of the boundary $ \, \partial \Omega $: 
\moneqstar  \left\{ \begin{array}{l} \displaystyle  
 \partial \Omega = \Gamma_m \cup \Gamma_p  \quad {\rm with}  \,\, 
{\rm meas} \, \big( \Gamma_m \cap \Gamma_p \big) = 0  \,, \\
\partial \Omega = \Gamma_t \cup \Gamma_\theta  \quad \,\, {\rm with}  \,\,  {\rm meas} \,
\big( \Gamma_t \cap \Gamma_\theta \big) = 0 \, . 
\end{array} \right.  \monendstar 
Suppose now that on one hand, 
the normal velocity is given on $ \, \Gamma_m \, $ and the pressure 
is  given on $ \, \Gamma_p \, $ and on the other hand that 
the tangential velocity   is given on $ \, \Gamma_t \, $ and the 
 tangential vorticity is imposed on  $ \, \Gamma_\theta $: 
\moneq   \label{stokes-cns-limites}  
\left\{ \begin{array}{ll} \displaystyle  
 u \cdot n = g_0  \quad {\rm on} \,\, \,   \Gamma_m \,,  \quad \,\, 
 & p = \Pi_0   \quad {\rm on} \,\, \,  \Gamma_p \\
 n \times u \times n = \sigma_0  \quad {\rm on} \,\, \,   \Gamma_t\,,  \quad \,\,  
 & n \times \omega  \times n = \theta_0  \quad {\rm on} \,\, \,  \, \Gamma_\theta  \,. 
\end{array} \right.  \monend   
Then under some  technical hypotheses, the Stokes problem  (\ref{stokes-omega-u-p}) 
with the boundary conditions (\ref{stokes-cns-limites}) admits a unique variational solution
in {\it ad hoc} vectorial Sobolev spaces \cite{Du02}. 
A first particular case is the Dirichlet problem, 
where both components of velocity are given (see {\it e.g.} \cite{GR86}).    
An  interesting case is the fact to give  pressure and tangential velocity, 
as first remarked in \cite{BCMP88}:
\moneq   \label{stokes-p-vtangenatielle}  
  p = \Pi_0   \quad {\rm and }  \quad 
 n \times u \times n = \sigma_0  \quad {\rm on} \,\, \,   \Gamma \, . 
\monend   
In the following of this section, we propose to adapt the algorithm proposed in Section~6
to the set of boundary conditions (\ref{stokes-p-vtangenatielle}). 
%
\begin{figure}    [H]  \centering 
  \centerline { \includegraphics[width=.35  \textwidth, angle=0] {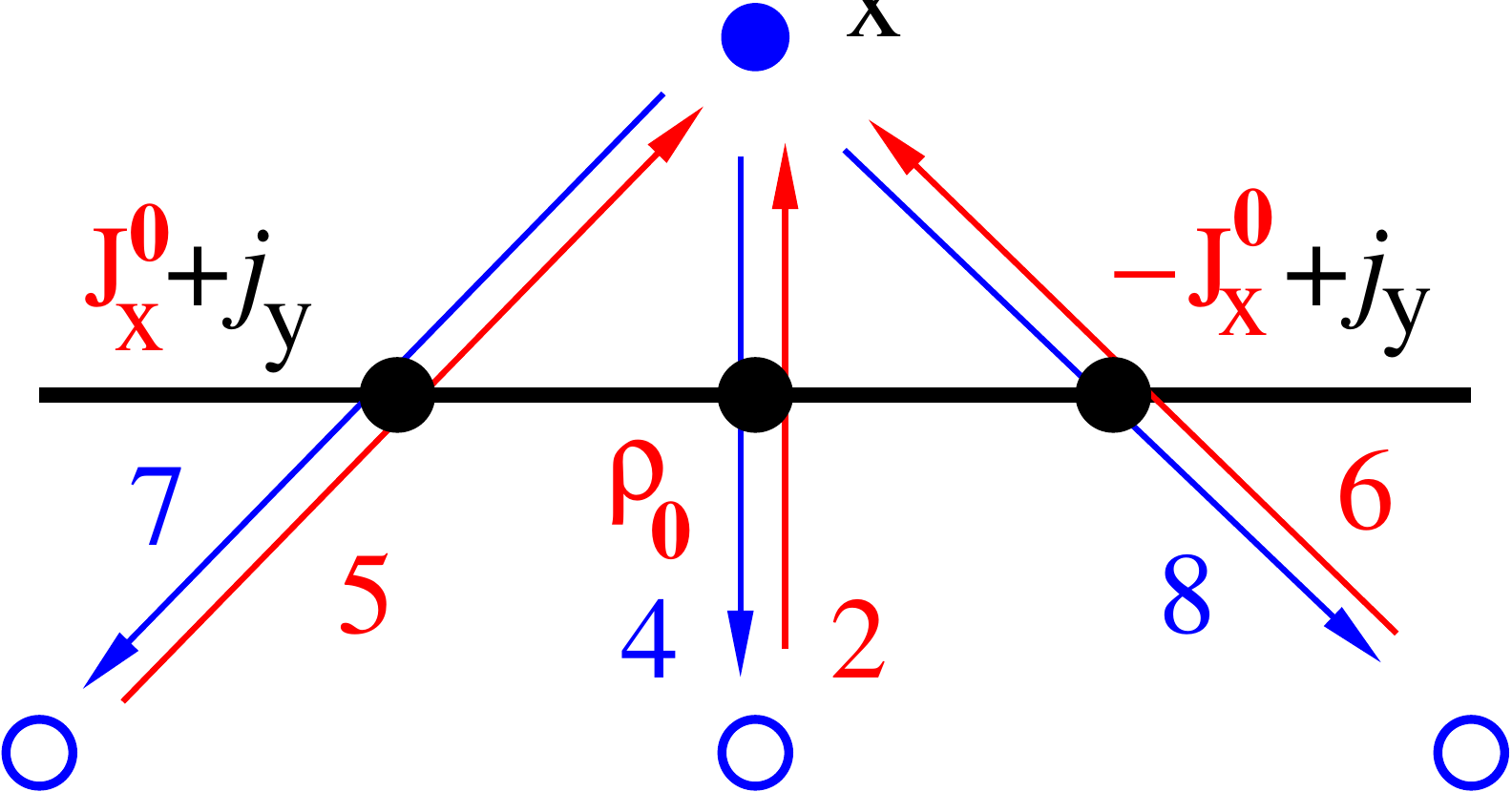}} 
\caption{Mixed bounce back and anti bounce back boundary condition to enforce 
density and tangential velocity on the boundary.    } 
\label{p-vtan-bb-abb} \end{figure}

\monitem 
We suppose in this section  that the density $ \, \rho = \rho_0 \, $ 
and the tangential momentum $ \, J_x = J_x^0 \, $ are given on the boundary. 
The momentum nodal values in the first cell are still denoted as $ \, j_x \, $ and $ \, j_y $. 
From the previous given equilibrium (\ref{feq-d2q9-fluide}), we can write 
the bounce back boundary condition for diagonal edges 
 $ \, f_5^{\rm eq} -  f_7^{\rm eq} = {{1}\over{6 \, \lambda}} \, ( J_x^0 + j_y ) $, 
$ \,   f_6^{\rm eq} -  f_8^{\rm eq} = {{1}\over{6 \, \lambda}} \, ( -J_x^0 + j_y ) $. 
The anti bounce back is unchanged: 
$ \,   f_2^{\rm eq} +  f_4^{\rm eq} = {{4 -\alpha-2 \, \beta}\over{18}} \, \rho_0 \, $ 
for the particles coming from the bottom, as illustrated in Figure~\ref{p-vtan-bb-abb}. 
If $ \, x \equiv \big(  x_1 , \, {{\Delta x}\over{2}} \big) \, $ is the notation 
for the first fluid node,    
the boundary condition on pressure and tangential velocity is finally implemented
in the D2Q9 algorithm with the relations 
\moneq \label{cl-p-vtan-bb-abb}  \left \{ \begin{array} {l} \displaystyle 
f_{5}(x, \, t+ \Delta t) =   f^*_{7}(x,\, t) \,+\,   {{1}\over{6 \, \lambda}} \, 
\Big[ J_x^0 \Big( x_1 - {{\Delta x}\over{2}}, \, t  \Big) +  j_y(x,\, t)   \Big]    \,, 
\\  \vspace{-4 mm} \\ \displaystyle
f_{2}(x,  \, t+ \Delta t) =    -    f^*_{4}(x,\, t) \,+\, {{4 -\alpha-2 \, \beta}\over{18}} \, 
\rho_0 \big(x_1 , \, t \big)  \,,  
\\  \vspace{-4 mm} \\ \displaystyle 
\displaystyle f_{6}(x,  \, t+ \Delta t)  =  f^*_{8}(x ,\, t) \,+\,  {{1}\over{6 \, \lambda}} \, 
\Big[ -J_x^0 \Big( x_1 + {{\Delta x}\over{2}}, \, t  \Big) +  j_y (x,\, t)   \Big]    \, . 
\end{array} \right.  \monend

\monitem 
The asymptotic analysis  of the  condition (\ref{cl-p-vtan-bb-abb})
follows the general framework ({\ref{cl-generale}). 
The matrix $ \, U \, $ is still given by  ({\ref{UU-cl-abb-chaleur}).
The interaction matrix $ \, T \, $ follows now  the relation  
\moneq \label{TT-cl-p-vtan-bb-abb} 
%
T = \left( \begin{array}{ccccccccc}
0 & 0 & 0 & 0 & 0 & 0 & 0 & 0 & 0 \\ 0 & 0 & 0 & 0 & 0 & 0 & 0 & 0 & 0 \\ 
0 & 0 & 0 & 0 & \!  -  1 & 0 & 0 & 0 & 0 \\ 0 & 0 & 0 & 0 & 0 & 0 & 0 & 0 & 0 \\ 
0 & 0 & 0 & 0 & 0 & 0 & 0 & 0 & 0 \\ 
0 & 0 & {1\over6}  & 0 &    \!  -{1\over6}  &   {1\over6}  &  {1\over6}  & \! {5\over6}  &  \!    -{1\over6} \\
0 & 0 &   {1\over6} & 0 &   \!  -{1\over6}  &  {1\over6}  &  {1\over6}  &   \!    -{1\over6}  &  \! {5\over6} \\ 
0 & 0 & 0 & 0 & 0 & 0 & 0 & 0 & 0  \\ 0 & 0 & 0 & 0 & 0 & 0 & 0 & 0 & 0 
\end{array}  \right) \, .  
\monend
The matrix $ \, K \, $ of  (\ref{matrice-K}) is now equal to   
\moneq \label{K-cl-p-vtan-bb-abb} 
K=\left( \begin{array}{ccccccccc} 
\!\!\!\! {{4-\alpha s_e-2 \beta s_d}\over{18}}  & \!\!\! \!\!\!\! 0 &  \!\!\!\! \!\!\!\!   {{1 - s_q }\over{3 \lambda}}
  &  \!\!\!   \!\!   {{s_e-1}\over{18  \lambda^2}} &  \!\!\!\!  {{s_x-1}\over{2  \lambda^2}} 
&  \!\!\!\! \!\!\! 0 &  \!\!\!  \!\!\!\! 0 &  \!\!\!\! {{1 - s_q}\over{3 \lambda^3}} 
&  \!\!\!\! {{s_d -1 }\over{9 \lambda^4}}   \\
\!\!\!\! 0 & \!\!\!  \!\!\!\!   {{2 - s_q }\over{3}}   & \!\!\!\!   \!\!\!\!  0 & \!\!\!   \!\!  0 &  \!\!\!\! 0 & 
 \!\!\!\! \!\!\!  0 &   \!\!\!\!  {{1 - s_q }\over{3 \lambda^2}}  &  \!\!\!\! 0 &  \!\!\!\! 0  \\ 
\!\!\!\!    {{\lambda  (4 - \alpha s_e - 2 \beta s_d )}\over{18}} & \!\!\! \!\!\!\!  0 &  
 \!\!\!\! \!\!\!\!    {{1 - s_q }\over{3}}   & 
 \!\!\!\!\!   {{s_e -1  }\over{18  \lambda}}  &  \!\!\!\!  {{s_x - 1}\over{2  \lambda}}  
&  \!\!\!\! \!\!\! 0 &  \!\!\!\! 0 &   \!\!\!\!  {{1 - s_q }\over{3 \lambda^2}}  &  
\!\!\!\!   {{s_d -1 }\over{9 \lambda^3}}   \\ 
\!\!\!\!  {{\lambda^2  (-4 - 17 \alpha s_e + 2 \beta s_d )}\over{18}} 
 & \!\!\! \!\!\!\!  0 &  \!\!\!\!  \!\!\!\!   {{2 \lambda (1 - s_q) }\over{3}} & \!\!\!   \!\!  {{1 + 17 s_e }\over{18}} 
 &  \!\!\!\!    {{1 - s_x}\over{2}} &  \!\!\! \!\!\!\! 0 &    \!\!\!\! 0 &  
\!\!\!\!  {{2 (1 - s_q) }\over{3 \lambda}}  &  \!\!\!\! {{1 - s_d }\over{9 \lambda^2}} \\ 
\!\!\!\!   {{\lambda^2  (-4 + \alpha s_e + 2 \beta s_d )}\over{18}} 
 & \!\!\! \!\!\!\! 0  &  \!\!\!\! \!\!\!\! 0  &  \!\!\!\!\!   {{1 - s_e }\over{18}} &  \!\!\!\!  {{1 + s_x}\over{2}} 
 &  \!\!\!\! \!\!\! 0  &  \!\!\!\! 0  &  \!\!\!\! 0 &  \!\!\!\! {{1 - s_d}\over{9 \lambda^2}}  \\ 
\!\!\!\! 0 & \!\!\!\!  \!\!\!\!  {{\lambda (2-s_q )}\over{3}} 
 &   \!\!\!\! \!\!\!\! 0 & \!\!\! \!\!   0 &  \!\!\!\! 0 & \!\!\!  \!\!\!\! s_x &
  \!\!\!\!   {{1-s_q }\over{3 \lambda}}   &  \!\!\!\!0 &  \!\!\!\!0 \\ 
\!\!\!\! 0 & \!\!\! \!\!\!\!  {{2 \lambda^2 (1+s_q) }\over3}   &  \!\!\!\!  \!\!\!\! 0 & \!\!\!  \!\!  0 &  \!\!\!\! 0 & 
\!\!\! \!\!\!\! 0  & \!\!\!\!   {{1 + 2 s_q  }\over3}  & \!\!\!\!0 & \!\!\!\!0 \\ 
\!\!\!\!   {{\lambda^3  (-4 + \alpha s_e + 2 \beta s_d )}\over{9}}  & \!\!\! \!\!\!\!  0 &  
 \!\!\!\! \!\!\!\!     {{\lambda^2 (1+ 2 s_q) }\over3} &
 \!\!\! \!\!  {{\lambda (1-s_e )}\over{9}} &  \!\!\!\! \lambda (1 \!\! - \!\! s_x\!) &  \!\!\! \!\!\!\! 0 & 
 \!\!\!\! 0 &  \!\!\!\!  {{1 + 2 s_q  }\over3} &  \!\!\!\!  {{2(1-s_d) }\over{9 \lambda}} \\ 
\!\!\!\!   {{\lambda^4  (-4 + \alpha s_e - 7 \beta s_d )}\over{9}}  &  \!\!\! \!\!\!\!  0 &  
 \!\!\!\! \!\!\!\!   {{\lambda^3 (1-s_q )}\over{3}}  &
\!\!\!  \!\! {{\lambda^2 (1-s_e)}\over{3}} &   \!\!\!\! \lambda^2 (1 \!\! - \!\! s_x \!) &  \!\!\!\! 0 &
 \!\!\!\! 0  & \!\!\! \!\!\!\!  {{\lambda (1-s_q)}\over{3}}  &  \!\!\!\! {{2 + 7 s_d}\over{9}} 
\end{array}  \right) \, . \monend

\smallskip 
 The matrix $K$ is  singular. Its   kernel   is of dimension 1. 
The null eigenvector of the matrix $ \, K \, $ is exactly $ \, \kappa_y \, $
introduced in (\ref{ker-K-fluide-abb}): 
\moneqstar
\kappa_y =  \big( 0 ,\, 0 ,\, 1 ,\, 0 ,\, 0 ,\, 0 ,\, 0 ,\, - \lambda^2 ,\, 0 \big)^{\rm  t}  \, . 
\monendstar 

\monitem 
The associated compatibility relation for solving the linear system (\ref{systeme-modele}) 
is analogous to previous ones:  
\moneq  \label{bbabb2-compatibilite} 
\lambda \, g_\rho - g_{jy} = 0  \, .
\monend
From the relation (\ref{cl-p-vtan-bb-abb}), the right hand side
\moneqstar
\xi = \xi_0 + \Delta t \,\, \partial \xi + {\rm O}(\Delta t^2) \, 
\monendstar
can be easily described:
%
%
\moneqstar \left\{ \begin{array}{ll} 
\displaystyle \xi_0  \! & = \displaystyle   \, \Big( 0 \,,\, 0 \,,\, {1\over18} \,(-\alpha - 2\, \beta + 4) \, \rho_0 \,,\, 0 \,,\, 0 \,,\,
{1\over{6 \, \lambda}} \,  J_x^0 \,,\, - {1\over{6 \, \lambda}} \,  J_x^0 \,,\,  0 \,,\, 0 \Big)^{\textrm t} \\ 
\partial \xi   \!  & = \displaystyle \, \Big( 0 \,,\, 0 \,,\, 0 \,,\, 0 \,,\, 0 \,,\,
- {1\over{12}} \,  \partial_x J_x^0 \,,\, - {1\over{12}} \,  \partial_x J_x^0 \,,\,  0 \,,\, 0 \Big)^{\! \textrm t} \, . 
\end{array} \right .  \monendstar 
At order zero, the compatibility relation (\ref{bbabb2-compatibilite}) 
does not   give any condition.   
At first order, the compatibility relation reduces to a combination 
of the associated partial differential  equations.
In consequence, no hidden boundary condition is introduced with this mixed bounce
back and anti bounce back boundary condition.
As designed by the boundary conditions (\ref{cl-p-vtan-bb-abb}), 
the normal momentum   $ \, j_y \, $    is   not defined  on the boundary.
Then all the results are defined up to a multiple $ \, j_y \, \kappa_y \, $ in the kernel
of the matrix $ \, K $. 
The  moments in the first cell at order zero are specified in the following relation: 
%
\moneqstar  
 m_0 =\left(\rho_0, \, J_x^0 , \,j_y , \,\alpha  \, \lambda^2\, \rho_0  ,\, 0,\, 0,\, 
- \lambda^2 \, J_x^0   ,\, - \lambda^2 \, j_y  , \,
\beta \, \,  \lambda^4 \rho_0  \right )^{\rm t}  \, . 
\monendstar
The boundary conditions $ \, \rho = \rho_0 \, $ and $ \, J_x  = J_x^0 \, $ are satisfied at  order zero
for the density and tangential momentum. 
At order one, the conserved moments $ \, \rho \, $ and $ \, j_x \, $ in the first cell have been expanded at order one.
Curiously, the  relations (\ref{bbabb-champs-premiere-maille}) obtained with the  approach in Section~7 are again valid
in this case:
\moneqstar 
\left\{ \begin{array}{l} \displaystyle 
  \rho \Big( x_1,\,  {{\Delta x}\over{2}} \Big)  = \rho_0 \big( x_1 \big) + {{1}\over{2}} \, \Delta x \, \partial_y \rho \big( x_1 ,\, 0 \big)
  +  {{ \Delta x}\over{s_x \, s_e \, s_d  \, (\alpha + 2\,\beta-4)}} \, \Big( a_\rho  \, \, \partial_y \rho \big( x_1 ,\, 0 \big)
\\ \vspace{-4 mm} \\ \displaystyle \qquad \qquad \qquad \qquad \qquad 
  + \, b_\rho \,\,   \partial_x J_x  \big( x_1 ,\, 0 \big) + c_\rho \,\,   \partial_y J_y  \big( x_1 ,\, 0 \big)  \Big)  +  {\rm O} (\Delta x^2) \,, 
\\ \vspace{-4 mm} \\ \displaystyle 
j_x  \Big( x_1, \, {{\Delta x}\over{2}} \Big) 
= J_x^0 \big( x_1 \big)   + {1\over2} \, \Delta x \, \, \partial_y J_x \big( x_1 ,\, 0 \big) +  {\rm O} (\Delta x^2) \,, 
\end{array} \right.  \monendstar 
with the coefficients $ \, a_\rho  \, $ and $ \, b_\rho \, $ given by the relations (\ref{bbabb-coefficients-champs}).

%
%

\monitem  Numerical experiments  for a linear Poiseuille flow 

\noindent
We have done the same numerical experiments than in Section 6 with a vertical linear Poiseuille flow.
The  pressure field is now converging towards the imposed linear profile
with a quasi first order accuracy, 
as described in Figure \ref{bbabb2-poiseuille-pression} and Table \ref{table-bbabb2-pression}.
We have observed  the horizontal component of the velocity.
A small discrepancy is present in the first layer (see Fig. \ref {bbabb2-poiseuille-vitesse-x}).
Nevertheless, the maximum error is reduced by one order of magnitude compared
  to the previous anti bounce back boundary  scheme. 
In the center of the channel, the error is very low: typically $ \, 1.5  \, 10^{-4} \, $
of relative error, that has to be compared to a relative error of  $ \, 10^{-4} \, $
with the previous pure anti bounce back.
The profile of vertical velocity (Fig. \ref{bbabb2-poiseuille-vitesse-y}) is quasi perfect.
The position of the numerical boundary, measured with the same technique  of least squares,
is presented in Table~\ref{table-bbabb}. These results are still not perfect
but can be considered  of good quality.

%
\begin{table}  [H]     \centering
 \centerline { \begin{tabular}{|c|c|c|c|}    \hline 
 mesh &  $ 20 \times 40 $ & $ 40 \times 80  $ & $   80 \times 160  $ \\   \hline 
 $ \, \nabla p $  & $  -1.06687  \, 10^{-3}   $ & $  -1.03617  \, 10^{-3} $ & $  -1.01899   \, 10^{-3} $   \\   \hline 
 relative error   &   6.69 \%    &   3.62 \%   &  1.90 \%    \\   \hline 
\end{tabular} }  
 \caption{Anti bounce back boundary condition for linear Poiseuille flow.
   Numerical gradient of the longitudinal pressure. The process is convergent with first order accuracy } 
\label{table-bbabb2-pression} \end{table}
%

\begin{figure}    [H]  \centering 
\vskip -1.8  cm
\centerline { \includegraphics[width=.95  \textwidth, angle=0] {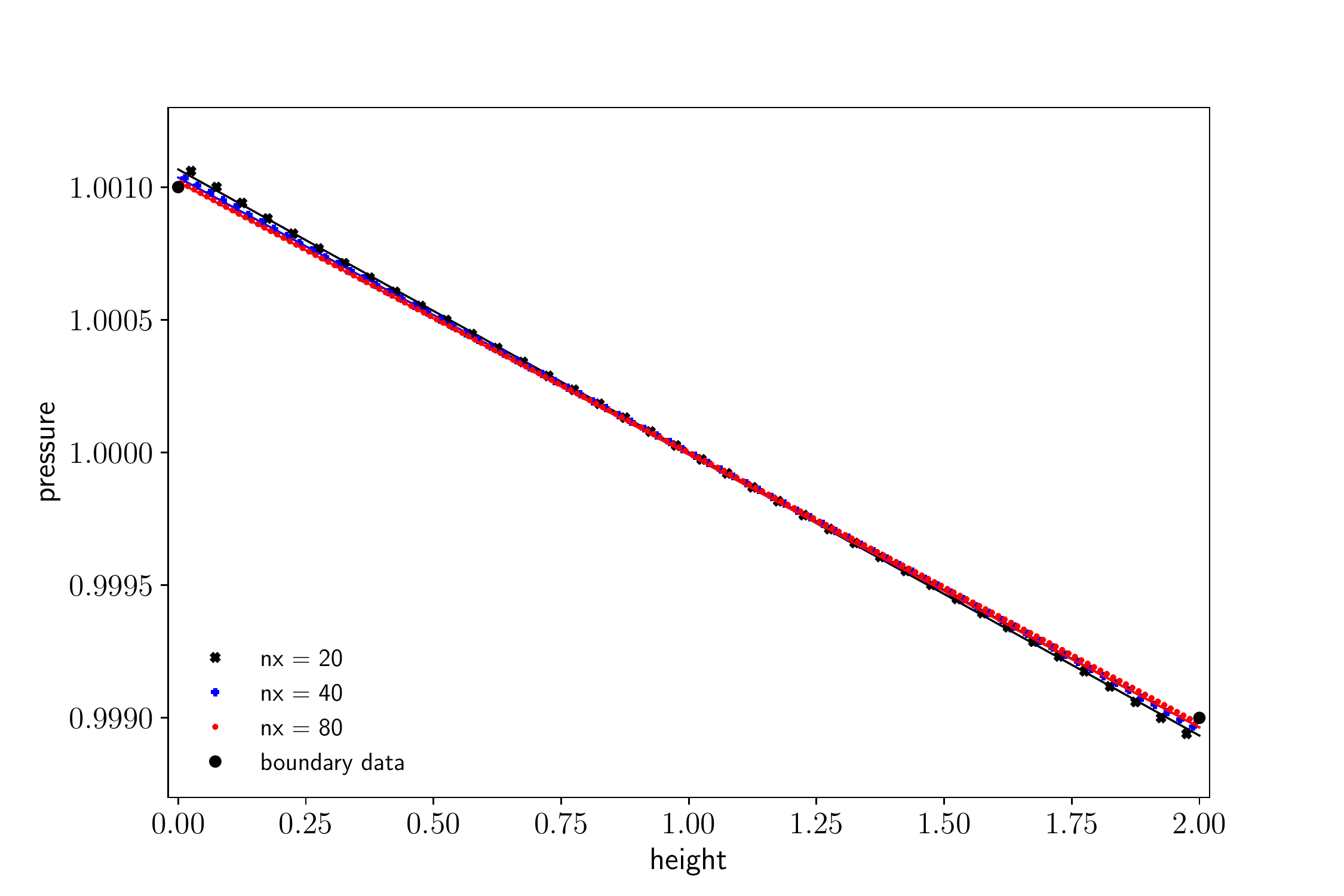}}
\vskip -.5 cm
\caption{Mixed bounce back and anti bounce back boundary condition for linear Poiseuille flow.
Pressure field for three meshes: the results are numerically convergent.   }
\label{bbabb2-poiseuille-pression} \end{figure}

\begin{figure}    [H]  \centering 
\vskip -.8  cm
\centerline { \includegraphics[width=.95  \textwidth, angle=0] {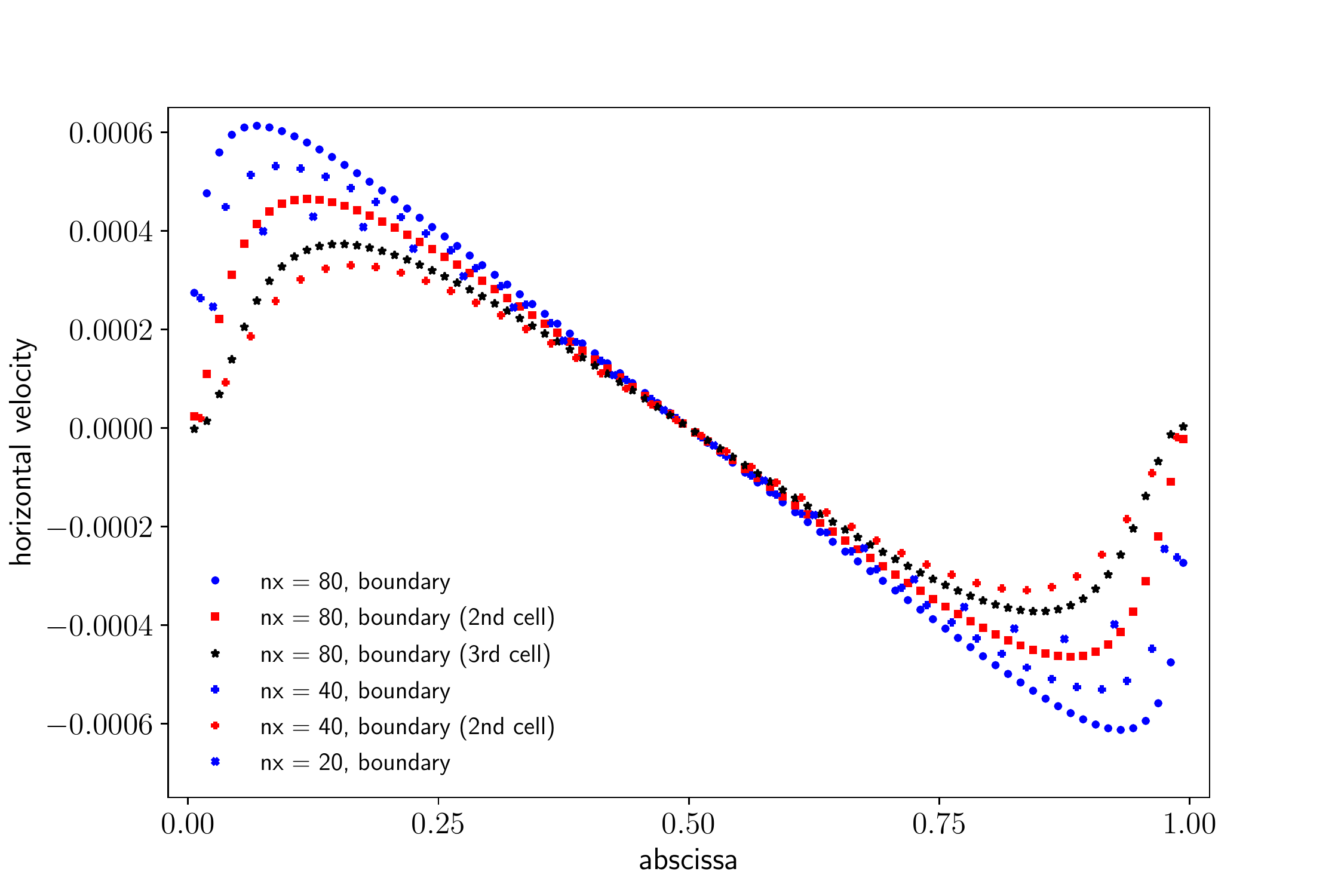}}
\vskip -.5 cm
\caption{Mixed bounce back and anti bounce back boundary condition for linear Poiseuille flow.
  Horizontal velocity field for three meshes. The maximum error is reduced by one order of magnitude compared
  to the anti bounce back boundary  scheme. } 
\label{bbabb2-poiseuille-vitesse-x} \end{figure}

\begin{figure}    [H]  \centering 
\vskip -1.8  cm
\centerline { \includegraphics[width=.95  \textwidth, angle=0] {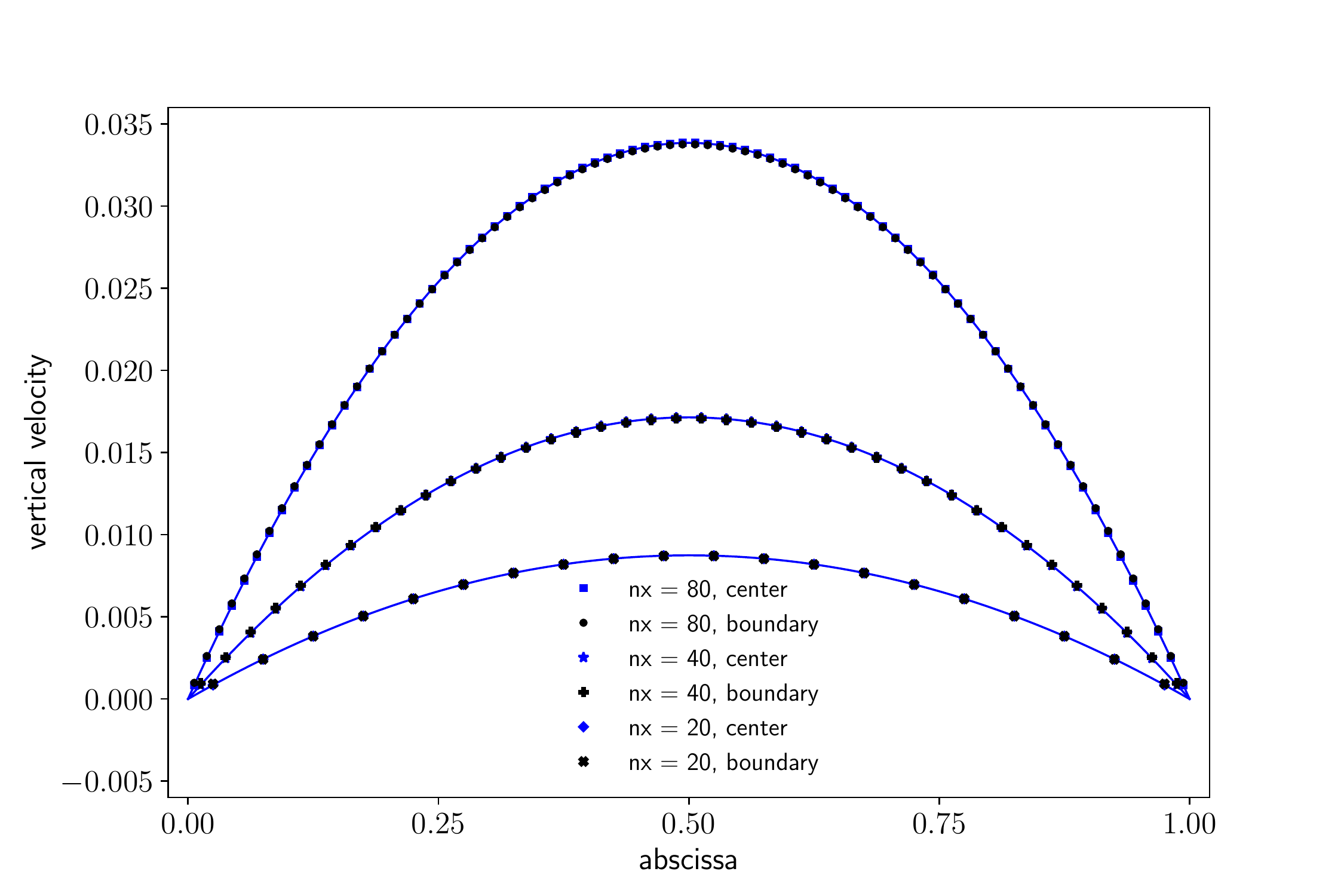}}
\vskip -.5  cm
\caption{Mixed bounce back and anti bounce back boundary condition for linear Poiseuille flow.
Vertical velocity field for three meshes. The parabolic profile is recovered in the first cell.   } 
\label{bbabb2-poiseuille-vitesse-y} \end{figure}

\vskip -.8  cm
%
\begin{table}  [H]     \centering
 \centerline { \begin{tabular}{|c|c|c|c|}    \hline 
 mesh &  $ 20 \times 40 $ & $ 40 \times 80  $ & $   80 \times 160  $ \\   \hline 
 center & $ 4.56  \, 10^{-3}   $ & $  2.40  \, 10^{-3}  $ & $  1.33  \, 10^{-3}  $  \\   \hline 
 bottom   & $  0.0165  $  & $  0.0608 $ &  $ 0.108 $  \\   \hline  
\end{tabular} }  
 \caption{Anti bounce back boundary condition for linear Poiseuille flow.
   Numerical position of the boundary measured in cell  units. The errors are reduced in a
   significative manner compared with the one obtained with the pure
 anti bounce boundary condition (see Table~\ref{table-abb}). } 
\label{table-bbabb} \end{table}
%

 \bigskip \bigskip    \noindent {\bf \large    9) \quad  Conclusion }   

In this contribution, we have presented the derivation of   anti bounce back boundary condition in the fundamental
case of linear heat conduction and linear acoustics. We have recalled that using this type of boundary
condition is now classical, even in the extended version that allows taking into account curved boundaries.
The asymptotic analysis confirms the high quality of the  anti bounce back boundary condition 
for  implementing a Dirichlet condition for the heat equation.
For the linear fluid system, the  anti bounce back boundary condition is designed for taking into
account a pressure boundary condition. 
The asymptotic analysis puts in evidence a  differential hidden condition on the boundary.
For a Poiseuille flow this hidden condition induces serious discrepancies in the vicinity of the
input and output. A variant mixing bounce back and anti  bounce back has been proposed in order to set in a mathematical rigorous
way fluid conditions composed by the datum of pressure and tangential velocity. A test case for Poiseuille flow
is very encouraging.

\noindent 
In future works, the analysis of the pure anti bounce back for thermal problem  will be extended  up to order two. 
The novel mixing boundary condition will be investigated more precisely theoretically  and numerically, 
in particular for unsteady acoustic problems.  
Moreover, the Taylor expansion method near the boundary needs to be generalized for
oblique and curved boundaries.  
Finally  a more general boundary condition has to be conceived to reduce  the defects at first order
and eliminate as far as possible the hidden boundary conditions.

%
\bigskip \bigskip   \noindent {\bf  \large  Acknowledgments }

\noindent 
The authors thank  the referees for precise comments on the first drafts of this contribution.


\bigskip \bigskip      \noindent {\bf  \large  References }   



\end{document}